\documentclass{article}
\usepackage[top=4cm, bottom=4cm, left=3cm, right=3cm]{geometry}
\usepackage{graphics}
\usepackage{epsfig}
\usepackage{enumerate}
\usepackage{amsmath, amssymb, amsfonts}
\usepackage{xcolor}
\usepackage[T1]{fontenc}
\usepackage{pifont}
\usepackage{subcaption}
\usepackage{graphicx,animate}
\usepackage{hyperref}
\usepackage{comment}

\numberwithin{equation}{section}

\newtheorem{theo}{Theorem}[section]

\newtheorem{prop}{Proposition}[section]
\newtheorem{lem}{Lemma}[section]
\newtheorem{definition}{Definition}[section]
\newtheorem{rem}{Remark}[section]

\newenvironment{proof}[1][Proof]{\textbf{#1.} }{\hfill \rule{0.5em}{0.5em}}

\def \R {\mathbb{R}}

\def \dis {\displaystyle}

\def \L {{\mathcal{L}}}
\def \U {{\mathcal{U}}}
\def \H {{\mathcal{H}}}

\begin{document}
  \title{Optimal Control of Covid-19 Interventions in Public Health Management}
\author{
	{{Isabella Kemajou-Brown} \thanks{{\it
				Department of Mathematics, EMERGE, NEERLab,  Morgan State University, Baltimore, MD 21251, USA,\\
				Email~: {\sf isabella.brown@morgan.edu} }}}
	{{~~ Romario Gildas Foko Tiomela} \thanks{{\it
				Department of Mathematics, EMERGE, Morgan State University, Baltimore, MD 21251, USA,\\
				Email~: {\sf romario.foko@morgan.edu} }}}
	{{~~ Olawale Nasiru Lawal} \thanks{{\it
				Department of Mathematics, EMERGE, Morgan State University, Baltimore, MD 21251, USA,\\
				Email~: {\sf ollaw5@morgan.edu} }}}\\
	{{~~ Samson Adekola Alagbe} \thanks{{\it
				Department of Mathematics, EMERGE, Morgan State University, Baltimore, MD 21251, USA,\\
				Email~: {\sf saala26@morgan.edu} }}}
	{{~~ Serges Love Teutu Talla} \thanks{{\it
				Department of Mathematics, EMERGE, Morgan State University, Baltimore, MD 21251, USA,\\
				Email~: {\sf seteu1@morgan.edu} }}}
} 
   
  \date{\today}
  
  \maketitle	
  	

\begin{abstract}
	This study explores the application of Pontryagin's Maximum Principle to derive optimal strategies for controlling the spread of COVID-19, leveraging a novel compartmental model to capture the disease dynamics. We prioritize three key criteria: cost, effectiveness, and feasibility, each examined independently to evaluate their unique contributions to pandemic management. By addressing these criteria, this study aims to design intervention strategies that are scientifically robust, practical, and economically sustainable. Furthermore, the focus on cost, effectiveness and feasibility seeks to provide policymakers with actionable insights for implementing interventions that maximize public health benefits while remaining feasible under real-world conditions.\\
\end{abstract}


\section{Introduction}
Optimal control theory plays a pivotal role in developing strategies to effectively manage infectious disease outbreaks. Among its foundational tools, Pontryagin's Maximum Principle provides a powerful framework for deriving optimal intervention strategies over time. The applications of Pontryagin's Maximum Principle span a wide range of contexts. H. W. Berhe \cite{berhe2018optimal} applied this principle to cholera outbreaks in Ethiopia, focusing on strategies for implementing treatment and sanitation to reduce the number of infected individuals and pathogenic agents while minimizing implementation costs. Similarly, M. A. Aba Oud et al. \cite{aba2021fractional} explored fractional-order mathematical models to assess the effects of quarantine, isolation, and environmental viral load on COVID-19 dynamics. Many other researchers have also investigated fractional optimal control problems in the context of disease dynamics, as discussed in \cite{baleanu2017new,baleanu2019new,jajarmi2019new,khan2020modeling}.\\

\noindent
Moreover, other researchers have further demonstrated the adaptability of Pontryagin's Maximum Principle in addressing various public health challenges. J. K. Asamoah et al. \cite{asamoah2020global} conducted sensitivity and cost-effectiveness analyses of a compartmental COVID-19 model, incorporating various control interventions. Extending this approach, J. K. Asamoah et al. \cite{asamoah2021sensitivity,asamoah2022optimal} emphasized the global stability and environmental impacts of control strategies, focusing on economic evaluations. L. Hakim \cite{hakim2023pontryagin} applied this principle to evaluate the effects of vaccination and quarantine strategies in reducing disease transmission. Additionally, K. O. Okosun et al. \cite{okosun2013optimal} analyzed malaria control strategies with an emphasis on cost-effectiveness, while S. Olaniyi et al. \cite{olaniyi2020mathematical} developed mathematical models to optimize COVID-19 interventions. Collectively, these works underscore the utility of Pontryagin's Maximum Principle in optimizing disease control efforts while balancing costs and outcomes. We refer to the following for additional resources \cite{abbasi2020optimal,lalwani2020predicting,mandal2020model,olaniyi2020modelling,tilahun2018co,ullah2020modeling,yousefpour2020optimal}.\\

\noindent
However, optimal control applications are not limited to single-disease contexts. For example, M. O. Adeniyi et al. \cite{adeniyi2018optimal} modeled the co-infection dynamics of malaria and pneumonia. Foundational insights into these methodologies are provided by L. Cesari \cite{cesari2012optimization} and L. S. Pontryagin \cite{pontryagin2018mathematical}, whose works have significantly advanced the theory of optimal control. E. B. Lee et al. \cite{lee1967foundations} also contributed to the theoretical foundations, offering key principles for deterministic systems, while A. D. Lewis \cite{lewis2006maximum} contextualized Pontryagin's Maximum Principle for both theoretical and applied settings.\\

\noindent
Beyond Pontryagin's Maximum Principle, data-driven approaches have emerged as important in epidemic modeling. For instance, M. Du et al. \cite{du2021data} proposed a multi-period resource allocation model for cholera control, using optimization to address the practical challenges of outbreak management. Machine learning and predictive analytics have further enriched this field. M. S. Ibrahim et al. \cite{ibrahim2023machine} highlighted their potential for advancing disease prevention, while O. V. Arueyingho et al. \cite{arueyingho2024scoping} reviewed their applications in healthcare systems management. Similarly, R. Yaesoubi et al. \cite{yaesoubi2016identifying} combined statistical and dynamic programming methods to identify cost-effective policies for controlling epidemics. In addition, M. U. Rahman et al. \cite{rahman2024mathematical} utilized computational fractional-order models to study the dynamics of communicable diseases, demonstrating their efficacy in capturing complex transmission patterns. Such approaches complement traditional methods by incorporating real-time data to enhance decision-making processes.\\

\noindent
Optimizing multiple criteria in the management of infectious disease outbreaks is critical for balancing health and economic outcomes. Particularly for complex and widespread illnesses like COVID-19, developing effective, sustainable, and cost-efficient public health strategies requires careful consideration of factors such as cost, feasibility, and intervention effectiveness. Each criterion plays a unique role in guiding decisions related to resource allocation, containment, and public cooperation. These factors shape a balanced strategy that not only minimizes the health impact of the virus but also considers broader societal and economic implications. This approach forms a foundation for informed decision-making in public health policy \cite{berhe2020optimal,fleming2012deterministic,kirk2004optimal,kuddus2024cost}.\\

\noindent
This article builds on previous research by applying Pontryagin's Maximum Principle to develop optimal strategies for controlling the spread of COVID-19. Unlike earlier models, this study explicitly incorporates three key criteria: cost, effectiveness, and feasibility, addressing the multifaceted challenges of pandemic management.\\

\noindent
The remainder of this article is structured as follows: Section \ref{prelim} introduces the mathematical preliminaries, including key concepts and theoretical results that form the foundation of this work. Section \ref{optim} details the proposed compartmental model, outlining its assumptions and structure, and develops the optimization framework using Pontryagin's Maximum Principle. Finally, Section \ref{conclusion} concludes this work and discusses its implications for public health decision-making.

\section{Preliminaries}\label{prelim}
In what follows, we consider the following spaces:
\begin{enumerate}[(i)]
	\item $C([a,b]):$ space of all continuous functions $h: [a,b] \longrightarrow \R$.
	\item $C^1([a,b]) = \left\{h\in C([a,b]): ~ h'\in C([a,b])\right\}$: space of all continuously differentiable functions $h: [a,b] \longrightarrow \R$.
	\item $PC([a,b]):$ space of all functions $h: [a,b] \longrightarrow \R$ that are piecewise continuous. In other words, $h(x)$ is continuous on each subinterval $[x_{i-1}, x_i]$ of a finite partition $a=x_0<x_1<\cdots<x_n=b$, and $h(x)$ possibly has a finite number of jumps.
	\item $PC^1([a,b]) = \left\{h\in PC([a,b]): ~ h|_{[x_{i-1}, x_i]} \in C^1([x_{i-1}, x_i]) ~\text{for all} ~ i \right\}:$ space of all functions $h:[a,b] \longrightarrow \R$ that are piecewise continuously differentiable. In other words, $h(x)$ is continuous on $[a,b]$, $h(x)$ is differentiable on each subinterval $[x_{i-1}, x_i]$, $h(x)$ is continuous on each subinterval $[x_{i-1}, x_i]$ and $h'(x)$ possibly has a finite number of jumps.
\end{enumerate}
\begin{definition}[Control System]
	A control system is a mathematical framework that describes how a dynamic process evolves under the influence of a control input. Formally, a control system $\Sigma=(\chi, f, \U)$ consists of the following components:
	\begin{enumerate}
		\item  A state space $\chi \subset \R^n$, which is an open subset of $\R^n$ containing all possible states of the system,
		\item A control set $\U\subset \R^m$, representing all admissible controls, and
		\item A function $f: \chi \times \overline{\U} \longrightarrow \R^n$, defining the system dynamics, and which is continuous and continuously differentiable (of class $C^1$) with respect to the state variable $x$ for each fixed control $u$ in the closed control set $\overline{\U}$.
	\end{enumerate}
	The system's evolution is captured by a differential equation:
	\begin{equation}
		\dis\frac{d x(t)}{d t} = f(x(t), u(t)), ~~ x(t_0)=x_0,
	\end{equation}
	where $u(t)\in \U$ is the control applied at time $t$. Control systems with these properties are foundational in optimal control theory and can be referenced in works such as E. B. Lee  et al. \cite{lee1967foundations} and D. E. Kirk \cite{kirk2004optimal}.
\end{definition}

\noindent
In this work, we will use the following important specific class of control system.

\begin{definition}[Control-affine system]
	A control-affine system is a particular class of control systems in which the state equation can be expressed as:
	\begin{equation}
		f(x,u) = f_0(x) + f_1(x)u,
	\end{equation}
	where $f_0: \chi \longrightarrow \R^n$ and $f_1: \chi \longrightarrow L(\R^m; \R^n)$ are of class $C^1$. Control-affine systems are commonly employed for systems where the control variables appear linearly, which is essential for utilizing Pontryagin's Maximum Principle.
\end{definition}

\noindent
A cornerstone of optimal control theory, Pontryagin's Maximum Principle \cite{pontryagin2018mathematical} (L. S. Pontryagin) provides first-order necessary conditions for optimality by transforming the control problem into a maximization or minimization problem for the Hamiltonian function.

\noindent
Let's consider the following control-affine problem:
\begin{equation}\label{PontMax}
	\begin{array}{rlllll}
		J(x(\cdot), u(\cdot)) &=& \dis\int_{a}^{b} g(t,x(t),u(t))dt + \phi(x(b)) \\~\\
		\text{subject to} ~ \dis\frac{dx(t)}{dt} &=& f(t,x(t),u(t)) ~ \text{and}~ x(a)=\alpha, ~ \alpha\in\R,
	\end{array}
\end{equation}
where $g\in C^1([a,b]\times\R^n\times \U; \R)$, the payoff term $\phi\in C^1(\R^n;\R)$, $f\in C^1([a,b]\times\R^n\times U; \R^n)$, the state variable $x\in PC^1([a,b];\R^n)$ and the control variable $u\in PC([a,b];\U)$, with $U\subseteq \R^r$ an open set.

\begin{prop}[Pontryagin's Maximum Principle]
	If a pair $(x(\cdot), u(\cdot))$ with $x\in PC^1([a,b];\R^n)$ and $u\in PC([a,b];\U)$ is a solution to problem \eqref{PontMax}, then there exists $\alpha\in PC^1([a,b];\R^n)$ such that the following conditions hold:
	\begin{enumerate}[(i)]
		\item the optimality condition
		\begin{equation}
			\dis\frac{\partial H}{\partial u}(t,x(t),u(t),\alpha(t)) = 0;
		\end{equation}
		\item the adjoint system
		\begin{equation}
			\left\{
			\begin{array}{rlllll}
				\dis\frac{dx(t)}{dt} &=& \dis\frac{\partial H}{\partial \alpha}(t,x(t),u(t),\alpha(t)) \\~\\
				\dis\frac{d\alpha(t)}{dt} &=& -\dis\frac{\partial H}{\partial x}(t,x(t),u(t),\alpha(t));
			\end{array}
			\right.
		\end{equation}
		\item and the transversality condition
		\begin{equation}
			\alpha(b) = \triangledown \phi(x(b));
		\end{equation}
	\end{enumerate}
	where the Hamiltonian $H$ is defined by
	\begin{equation}
		H(t,x,u,\alpha) = g(t,x,u) + \alpha\cdot f(t,x,u)
	\end{equation}
	and $\alpha\in\R^n$ is the adjoint variable (or costate) associated with the state $x$.
\end{prop}

\noindent
\begin{rem}
	This principle asserts that an optimal control $u^*$ should satisfy:
	\begin{equation}
		H(t,x^*(t),u^*(t),\alpha^*(t)) = \max_{u\in\U} H(t,x^*(t),u,\alpha^*(t)),
	\end{equation}
	or
	\begin{equation}
		H(t,x^*(t),u^*(t),\alpha^*(t)) = \min_{u\in\U} H(t,x^*(t),u,\alpha^*(t)),
	\end{equation}
	where $x^*(t)$ and $\alpha^*(t)$ are the optimal state and adjoint trajectories.
\end{rem}
We refer to L. S. Pontryagin \cite{pontryagin2018mathematical} and A. D. Lewis \cite{lewis2006maximum} for more.

\begin{definition}[Lagrange function]
	The Lagrange function $L(t,x,u)$ (or Lagrangian) is a term in the objective functional that represents the instantaneous cost or reward of being in a particular state $x$ while applying control $u$ at time $t$. The objective functional is given by:
	\begin{equation}
		J(x(\cdot), u(\cdot)) = \dis\int_{a}^{b} L(t,x(t),u(t))dt,
	\end{equation}
	where $L: [a,b]\times \R^n \times \U$ is the Lagrangian. For more information on Lagrange functions, we refer to W. H. Fleming et al. \cite{fleming2012deterministic}.
\end{definition}

\noindent
Several important results ensure the existence of optimal control solutions for systems like the one studied in this article.

\begin{prop}[Existence of an optimal control]\label{ExistContr}
	Let $\Sigma = (\chi, f, \U)$ be a control system where $\chi \subset \R^n$ is the state space, $\U\subset \R^m$ is a closed and convex control set, and $f: \chi \times \U \longrightarrow \R^n$ is a continuously differentiable function, representing the system dynamics and governed by the differential equation:
	\begin{equation}
		\dis\frac{dx}{dt} = f(x(t),u(t)), ~~ x(0)=x_0.
	\end{equation}
	
	\noindent
	Consider an objective functional $J$ defined by:
	\begin{equation}
		J(u) = \dis\int_{0}^{T} L(t,x(t), u(t))dt + \phi(x(T)),
	\end{equation}
	where $L: [0,T]\times \R^n\times\R^m \longrightarrow \R$ is the Lagrange function, and $\phi: \R^n \longrightarrow \R$ is the terminal cost function, both of which are continuously differentiable. If the following conditions hold:
	\begin{enumerate}[(i)]
		\item $\U$ is closed, convex, and bounded;
		\item $f(x,u)$ is linearly bounded in $x$ and $u$;
		\item $L(t,x,u)$ is convex (resp. concave) in $u$ for each fixed $(t,x)$;
		\item $L$ is bounded below (resp. above) by $L(t,x,u) \geq a_0\left(\dis\sum_{i=1}^{m}u_i^2\right)^{\frac{a_2}{2}}-a_1$ (resp. $L(t,x,u) \leq b_0\left(\dis\sum_{i=1}^{m}u_i^2\right)^{\frac{b_2}{2}}+b_1$) for constants $a_0, b_0, a_1, b_1>0$ and $a_2, b_2>1$,
	\end{enumerate}
	 then there exists an optimal control $u^*\in \U$ such that:
	 \begin{equation}
	 	J(u^*) = \min_{u\in\U} J(u)
	 \end{equation}
	resp.
	 \begin{equation}
	 	J(u^*) = \max_{u\in\U} J(u).
	 \end{equation}
\end{prop}
For a detailed discussion and proof of this result, we refer to W. H. Fleming et al. \cite{fleming2012deterministic}, and L. Cesari \cite{cesari2012optimization}

\newpage
\section{Optimal control of a COVID-19 model}\label{optim}

\subsection{Model Description}\label{SectDescript}
\begin{figure}[h!]
	\centering
	\includegraphics[height=6cm]{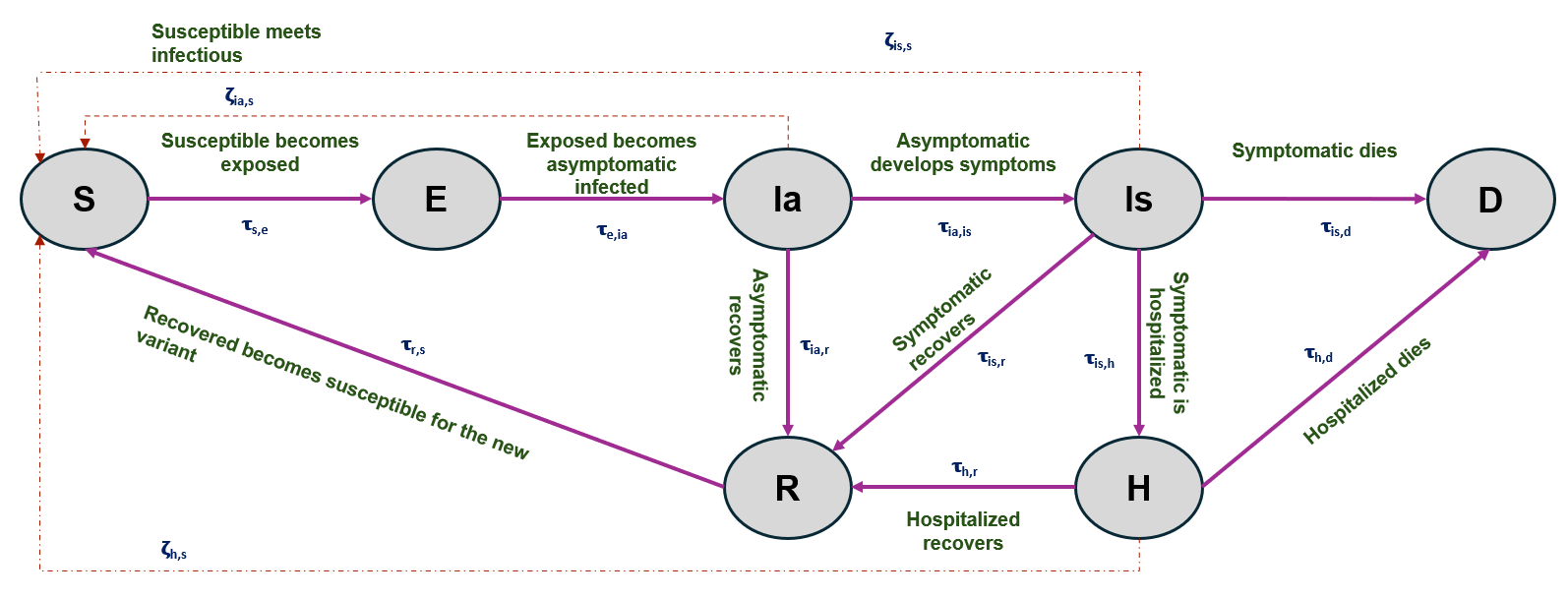}
	\caption{Dynamics of COVID-19 spread}
	\label{FigModel}
\end{figure}

Our model (Figure \ref{FigModel}) simulates the spread of COVID-19 through a population, distributed into different compartments, each representing a distinct stage of the disease.

\noindent
This model is typically represented using the following system of differential equations:

\begin{equation}\label{diff}
	\left\{
	\begin{array}{rlllll}
		\dis\frac{dS}{dt}&=& -\left(\zeta_{ia,s}I_a + \zeta_{is,s}I_s + \zeta_{h,s}H\right)\dis\frac{S}{N} + \tau_{r,s}R \\~\\
		\dis\frac{dE}{dt}&=& \left(\zeta_{ia,s}I_a + \zeta_{is,s}I_s + \zeta_{h,s}H\right)\dis\frac{S}{N} - \tau_{e,ia}E \\~\\
		\dis\frac{dI_a}{dt}&=& \tau_{e,ia}E - (\tau_{ia,is} + \tau_{ia,r})I_a \\~\\
		\dis\frac{dI_s}{dt}&=& \tau_{ia,is}I_a - (\tau_{is,r} + \tau_{is,h} + \tau_{is,d})I_s \\~\\
		\dis\frac{dH}{dt}&=& \tau_{is,h}I_s - (\tau_{h,r} + \tau_{h,d})H \\~\\
		\dis\frac{dR}{dt}&=& \tau_{is,r}I_s + \tau_{ia,r}I_a + \tau_{h,r}H - \tau_{r,s}R \\~\\
		\dis\frac{dD}{dt}&=& \tau_{h,d}H + \tau_{is,d}I_s
	\end{array}
	\right.,
\end{equation}
where $S, E, I_a, I_s, H, R$ and $D$ represent the states of Susceptible, Exposed, Asymptomatic Infected, Symptomatic Infected, Hospitalized, Recovered and Deceased individuals respectively. The different states are described as follows:

\begin{enumerate}[(i)]
	\item Susceptible: Individuals who have not been infected with the virus but are at risk of becoming infected. They can move to the Exposed compartment upon contact with an infected individual.
	
	\item Exposed: Individuals who have been exposed to the virus and are in the incubation period. They are infected but not yet infectious. After the incubation period, they transition to the Asymptomatic Infected compartment.
	
	\item Asymptomatic Infected: Individuals who have been infected with the virus but do not show symptoms. They can still transmit the virus to Susceptible individuals. After a certain period, they either recover or transition to the symptomatic infected compartment.
	
	\item Symptomatic Infected: Individuals who show symptoms of the infection. These individuals can be further categorized by the severity of symptoms. Depending on the severity and progression of the disease, they may either recover, be hospitalized, or eventually transition to the deceased compartment if the infection becomes fatal.
	
	\item Hospitalized: Individuals with severe symptoms who require intensive medical care. Their outcomes can vary: they may recover and move to the Recovered compartment or, in severe cases, may not survive and move to the Deceased compartment.
	
	\item Recovered: Individuals who have recovered from the infection but their immunity can last only for some period. They cannot transmit the disease but are still at risk of infection.
	
	\item Deceased: Individuals who have died due to the infection.
\end{enumerate}

\noindent
The transitions between different compartments of the model \eqref{diff} are governed by probabilities, which can be stationary (constant over time) or non-stationary (varying over time). These transition probabilities depend on various factors, including the disease's natural progression, intervention measures, and population behavior, and were comprehensively discussed in a previous study. Our state transitions are described as follows:
\begin{enumerate}[(a)]
	\item $\tau_{s,e}:$ probability for a susceptible individual to be exposed to the disease. It depends on the contact rate with infected individuals (asymptomatic, symptomatic or hospitalized), influenced by factors like social distancing, mask usage, vaccination rates, and so on.
	
	\item $\tau_{e,ia}:$ probability for an individual to be infected by the disease without exhibiting symptoms. This can be influenced by demographic factors such as underlying health conditions.
	
	\item $\tau_{ia,r}:$ probability for an individual to recover from asymptomatic infection. This probability is more likely based on the average duration of the infection.
	
	\item $\tau_{ia,is}:$ probability for an individual to develop symptoms subsequent to previous asymptomatic infection.
	
	\item Symptomatic individuals can have different pathways: mild cases may recover without hospitalization, severe cases may require hospitalization, followed by recovery or death, fatal cases transition to deceased.
	
	\begin{itemize}
		\item $\tau_{is,h}:$ probability for an individual to require hospitalization in the Intensive Care Unit (ICU) due to complications arising from the symptoms of the disease.
		
		\item $\tau_{is,r}:$ probability for an individual to successfully recover after exhibiting symptoms of the disease.
		
		\item $\tau_{is,d}:$ probability for an individual to experience complications leading to death after displaying symptoms of the disease.
	\end{itemize}
	
	\item Hospitalized individuals may recover or succumb to the disease, with probabilities influenced by the quality of healthcare, the severity of the disease, and comorbidities.
	
	\begin{itemize}
		\item $\tau_{h,d}:$ probability for an individual to succumb to mortality following hospitalization in the ICU.
		
		\item $\tau_{h,r}:$ probability for an individual to successfully recover after being discharged from the ICU.
	\end{itemize}
	
	\item $\tau_{r,s}:$ probability for an individual to become susceptible again after recovering from the disease.
	
	\item $\zeta_{ia,s}:$ probability to contract the infection in one meeting with an individual from $I_a$.
	
	\item $\zeta_{is,s}:$ probability to contract the infection in one meeting with an individual from $I_s$.
	
	\item $\zeta_{h,s}:$ probability to contract the infection in one meeting with an individual from $H$.
\end{enumerate}

\begin{rem}
	We note that susceptible individuals become exposed by coming into contact with asymptomatically infectious, symptomatically infectious and hospitalized individuals.
\end{rem}

\noindent
To act on system \eqref{diff}, we introduce the time-dependent control vector $u = (u_1, u_2, \cdots, u_9)$ described in the table below:

\begin{table}[h!]
	\centering
	\caption{Control Description}
	\begin{tabular}{|c|c|}
	\hline
	control variable &  description\\
	\hline
	$u_1$ &  vaccination \\
	\hline
	$u_2$ &  nothing \\
	\hline
	$u_3$ &  social distancing \\
	\hline
	$u_4$ &  hand sanitization \\
	\hline
	$u_5$ &  masking \\
	\hline
	$u_6$ &  testing and contact tracing of exposed individuals \\
	\hline
	$u_7$ &  admission to the hospital \\
	\hline
	$u_8$ &  quarantine \\
	\hline
	$u_9$ &  Intensive medical care tailored for individuals receiving treatment in the hospital \\
	\hline
	\end{tabular}\label{Contr_var}
\end{table}

\noindent
The corresponding control system describes how the different compartments evolve over time, as follows:

\begin{equation}\label{contsyst1}
	\left\{
	\begin{array}{rlllll}
		\dis\frac{dS}{dt} & = & -\left(1-u_1-u_3-u_4-u_5\right) \left(\zeta_{ia,s}I_a + \zeta_{is,s}I_s + \zeta_{h,s}H\right)\dis\frac{S}{N} \\
		&& - u_2\left(\zeta_{ia,s}I_a + \zeta_{is,s}I_s + \zeta_{h,s}H\right)\dis\frac{S}{N} + (1-u_3-u_4-u_5)\tau_{r,s}R \\~\\
		\dis\frac{dE}{dt} & = & \left(1-u_1-u_3-u_4-u_5\right) \left(\zeta_{ia,s}I_a + \zeta_{is,s}I_s + \zeta_{h,s}H\right)\dis\frac{S}{N}\\
		&& + u_2\left(\zeta_{ia,s}I_a + \zeta_{is,s}I_s + \zeta_{h,s}H\right)\dis\frac{S}{N} - (1-u_6)\tau_{e,ia}E \\~\\
		\dis\frac{dI_a}{dt} & = & (1-u_6)\tau_{e,ia}E - (1-u_6) (\tau_{ia,is} + \tau_{ia,r})I_a \\~\\
		\dis\frac{dI_s}{dt} & = & (1-u_6)\tau_{ia,is}I_a - (1-u_7-u_8)(\tau_{is,r} + \tau_{is,h} + \tau_{is,d})I_s \\~\\
		\dis\frac{dH}{dt} & = & (1-u_7-u_8)\tau_{is,h}I_s - (1-u_9)\tau_{h,d}H - u_9\tau_{h,r}H \\~\\
		\dis\frac{dR}{dt} & = & (1-u_7-u_8)\tau_{is,r}I_s + (1-u_6)\tau_{ia,r}I_a + u_9\tau_{h,r}H \\
		&& -(1-u_3-u_4-u_5)\tau_{r,s}R \\~\\
		\dis\frac{dD}{dt} & = & (1-u_9)\tau_{h,d}H + (1-u_7-u_8)\tau_{is,d}I_s,
	\end{array}
	\right.
\end{equation}
where the control variables $u_1, u_2, \cdots, u_9$ represent the various intervention strategies described in Table \ref{Contr_var}.

\noindent
We assume that the control set $\U$ is defined by:

$$\U = \left\{u= (u_1, \cdots, u_9): ~ u_i\big|_{i=1}^{9} ~ \text{is Lebesgue measurable}, ~~ 0\leq u_i \leq 1, ~ i=1, \cdots, 9, ~ t\in [0,T] \right\},$$
where the time horizon $T$ is fixed, and there is no explicit dependency on the state variable at $T$.
If we denote 
$$\max_{i=1,\cdots,9} \left\{u_i\right\} = u_{max},$$
then
$$0 \leq u_i \leq u_{max}$$
which means that $u_i\in [0, u_{max}]$, $i=1,2, \cdots, 9$. It then follows that
$$u = (u_1, \cdots, u_9) \in [0, u_{max}]^9.$$

\noindent
In disease control, understanding the impact of cost, effectiveness and feasibility on an intervention strategy is crucial for informed decision-making. Optimizing each of these criteria offers valuable insights into the distinct contributions and limitations of each criterion. This approach allows for a focused analysis on how minimizing or maximizing each criterion impacts the overall outcome of COVID-19 control strategy. This way, it becomes possible to isolate and fully understand the effects that each criterion has on the system. For example, when cost is minimized independently, we can identify the lowest possible financial commitment required to control the outbreak. Similarly, by maximizing effectiveness alone, we can understand the maximum achievable impact in terms of reducing infection and mortality rates, independent of resource or logistical constraints. Maximizing feasibility by itself can clarify which interventions are most realistic and achievable within the limitations of current healthcare infrastructure and public compliance. In what follows, we optimize each criterion individually.

\subsection{Cost Objective Function}
The cost criterion represents the financial and resource expenditures associated with implementing a particular set of control measures. For large-scale outbreaks, these costs can be considerable, and resource constraints often limit the extent of measures that can be feasibly deployed. Optimizing the cost criterion involves minimizing expenditures while ensuring that the implemented measures are effective enough to contain the spread of the disease. This approach enables policymakers to allocate limited resources wisely, potentially extending the reach of interventions to more communities and allowing for sustained efforts over time. Let's consider the cost objective function
\begin{equation}\label{objective1}
	J_1(u_1, u_2, \cdots, u_9) = \dis\int_{0}^{T} \left[ \lambda_1E + \lambda_2I_a + \lambda_3I_s + \lambda_4H + \dis\frac{1}{2}\dis\sum_{i=1}^{9}b_iu_i^2(t)\right]dt,
\end{equation}
where $\left\{\lambda_i\right\}_{i=1}^{4}$ are the balancing weight on the exposed, the asymptomatically infected, the symptomatically infected and the hospitalized individuals; $\left\{b_i\right\}_{i=1}^{9}$ are the balancing cost factors on the different controls $\left\{u_i\right\}_{i=1}^{9}$; and $T$ is the final time for control implementation.

\noindent
	Our goal is to \textbf{find an optimal control} $u^* = (u_1^*, u_2^*, \cdots, u_9^*)$ such that 
	\begin{equation}
		J_1(u^*) = \min_{u\in\U}{J_1(u)}
	\end{equation}
	where $u=(u_1,u_2, \cdots, u_9) \in \U$.

\noindent
Let $x=(S,E,I_a,I_s,H,R,D)$ be a solution to the state system \eqref{contsyst1}, and consider
\begin{equation}\label{affine-funct}
	f_0(t,x,u) = 
	\begin{pmatrix}
		-\left(1-u_1+u_2-u_3-u_4-u_5\right) \left(\zeta_{ia,s}I_a + \zeta_{is,s}I_s + \zeta_{h,s}H\right)\dis\frac{S}{N} \\
		+ (1-u_3-u_4-u_5)\tau_{r,s}R \\~\\
		\left(1-u_1+u_2-u_3-u_4-u_5\right) \left(\zeta_{ia,s}I_a + \zeta_{is,s}I_s + \zeta_{h,s}H\right)\dis\frac{S}{N} - (1-u_6)\tau_{e,ia}E \\~\\
		(1-u_6)\tau_{e,ia}E - (1-u_6) (\tau_{ia,is} + \tau_{ia,r})I_a \\~\\
		(1-u_6)\tau_{ia,is}I_a - (1-u_7-u_8)(\tau_{is,r} + \tau_{is,h} + \tau_{is,d})I_s \\~\\
		(1-u_7-u_8)\tau_{is,h}I_s - (1-u_9)\tau_{h,d}H - u_9\tau_{h,r}H\\~\\
		(1-u_7-u_8)\tau_{is,r}I_s + (1-u_6)\tau_{ia,r}I_a + u_9\tau_{h,r}H \\
		-(1-u_3-u_4-u_5)\tau_{r,s}R \\~\\
		(1-u_9)\tau_{h,d}H + (1-u_7-u_8)\tau_{is,d}I_s
	\end{pmatrix}
\end{equation}

\noindent
The vector $f_0$ can be decomposed by $f_0(t,x,u) = f_1(t,x) + f_2(t,x)u$, where:

\begin{equation*}
	f_1(t,x) = 
	\begin{pmatrix}
		-\left(\zeta_{ia,s}I_a + \zeta_{is,s}I_s + \zeta_{h,s}H\right)\dis\frac{S}{N} + \tau_{r,s}R \\~\\
		\left(\zeta_{ia,s}I_a + \zeta_{is,s}I_s + \zeta_{h,s}H\right)\dis\frac{S}{N} - \tau_{e,ia}E \\~\\
		\tau_{e,ia}E - (\tau_{ia,is} + \tau_{ia,r})I_a \\~\\
		\tau_{ia,is}I_a - (\tau_{is,r} + \tau_{is,h} + \tau_{is,d})I_s \\~\\
		\tau_{is,h}I_s - \tau_{h,d}H \\~\\
		\tau_{is,r}I_s + \tau_{ia,r}I_a - \tau_{r,s}R \\~\\
		\tau_{h,d}H + \tau_{is,d}I_s
	\end{pmatrix}
\end{equation*}
and
\begin{equation*}
	f_2(t,x) =
	\begin{pmatrix}
		c_{11} & c_{12} & \cdots & c_{19} \\~\\
		c_{21} & c_{22} & \cdots & c_{29} \\~\\
		\vdots & \vdots & \ddots & \vdots \\~\\
		c_{71} & c_{72} & \cdots & c_{79}
	\end{pmatrix}
\end{equation*}
where:
$c_{11} = (\zeta_{ia,s}I_a + \zeta_{is,s}I_s + \zeta_{h,s}H)\dis\frac{S}{N}$, $c_{12} = -c_{11}$, $c_{13} = c_{11} - \tau_{r,s}R$, $c_{14} = c_{15} = c_{13}$, $c_{16} = c_{17} = c_{18} = c_{19} = 0$, $c_{21} = c_{23} = c_{24} = c_{25} = -c_{11}$,  $c_{22} = c_{11}$, $c_{26} = \tau_{e,ia}E$, $c_{27} = c_{28} = c_{29} = 0$, $c_{31} = c_{32} = c_{33} = c_{34} = c_{35} = c_{37} = c_{38} = c_{39} = 0$, $c_{36} = -\tau_{e,ia}E + (\tau_{ia,is}+\tau_{ia,r})I_a$, $c_{41} = c_{42} = c_{43} = c_{44} = c_{45} = c_{49} = 0$, $c_{46} = -\tau_{ia,is}I_a$, $c_{47} = c_{48} = (\tau_{is,r} + \tau_{is,h} + \tau_{is,d})I_s$, $c_{51} = c_{52} = c_{53} = c_{54} = c_{55} = c_{56} = 0$, $c_{57} = c_{58} = -\tau_{is,h}I_s$, $c_{59} = (\tau_{h,d}-\tau_{h,r})H$, $c_{61} = c_{62} = 0$, $c_{63} = c_{64} = c_{65} = \tau_{r,s}R$, $c_{66} = -\tau_{ia,r}I_a$, $c_{67} = c_{68} = - \tau_{is,r}I_s$, $c_{69} = \tau_{h,r}H$, $c_{71} = c_{72} = c_{73} = c_{74} = c_{75} = c_{76} = 0$, $c_{77} = c_{78} = -\tau_{is,d}I_s$, $c_{79} = -\tau_{h,d}H$.\\

\noindent
Now, we can state and prove the following existence result of an optimal control.

\begin{theo}\label{control1}
	Given the objective function $J_1$ defined in \eqref{objective1} on the control set $\U$, there exists an optimal control $u^* = (u_1^*, u_2^*, \cdots, u_9^*)$ such that 
	\begin{equation}\label{optimiz1}
		J_1(u^*) = \min_{u\in\U}{J_1(u)}
	\end{equation}
	provided that the following conditions are all satisfied.
	
	\begin{enumerate}[(i)]
		\item The admissible control set $\U$ is closed and convex.
		\item The function $f_0(t,x,u)$ defined in \eqref{affine-funct} is bounded by a linear function in the state and control variables.
		\item The integrand 
		$$\L_1(t,x,u) = \lambda_1 E + \lambda_2 I_a + \lambda_3 I_s + \lambda_4 H + \dis\frac{1}{2}\dis\sum_{i=1}^{9}b_iu_i^2,$$
		of the objective function $J_1$ is convex with respect to the control $u$.
		\item The Lagrangian $\L_1(t,x,u)$ is bounded below by $\dis\frac{1}{2}\dis\min_{i=1,\cdots,9}\left\{b_i\right\}\dis\sum_{i=1}^{9}u_i^2 - a_1$
		where $a_1$ is any positive constant.
	\end{enumerate}
\end{theo}

\noindent
\begin{proof}
\begin{itemize}
	\item[(i)] \textbf{Prove that the admissible control set $\U$ is closed and convex.}\\
	By construction, it is clear that $\U \subset [0, u_{max}]^9$; and hence, \textbf{closed} and \textbf{convex}.
	\item[(ii)] \textbf{Prove that $f_0(t,x,u)$ is bounded by a linear function in the state and control variables.}\\
	
	We know that $f_0(t,x,u) = f_1(t,x) + f_2(t,x)u$. By Triangle Inequality,
	$$\|f_0(t,x,u)\| \leq \|f_1(t,x)\| + \|f_2(t,x)\| \|u\|.$$
	
	Now, we have:
	
\end{itemize}

$$\begin{array}{rlll}
	&&\|f_1(t,x)\|^2 \\
	
	&\leq& 2(\zeta_{ia,s}I_a + \zeta_{is,s}I_s + \zeta_{h,s}H)^2\dis\frac{S^2}{N^2} + \tau_{r,s}^2R^2 - 2\tau_{r,s}(\zeta_{ia,s}I_a + \zeta_{is,s}I_s + \zeta_{h,s}H)\dis\frac{S}{N} \\~\\
	
	&& + \tau_{e,ia}^2E^2 - 2\tau_{e,ia}(\zeta_{ia,s}I_a + \zeta_{is,s}I_s + \zeta_{h,s}H)\dis\frac{S}{N} + \tau_{e,ia}^2E^2 + (\tau_{ia,is}+\tau_{ia,r})^2I_a^2\\~\\
	
	&& - 2\tau_{e,ia}(\tau_{ia,is}+\tau_{ia,r})EI_a + \tau_{ia,is}^2I_a^2 + (\tau_{is,r}+\tau_{is,h}+\tau_{is,d})^2I_s^2 \\~\\
	
	&&-2\tau_{ia,is}(\tau_{is,r}+\tau_{is,h}+\tau_{is,d})I_aI_s +\tau_{is,h}^2I_s^2 + \tau_{h,d}^2H^2 \\~\\
	
	&& -2\tau_{is,h}\tau_{h,d}I_sH + \tau_{is,r}^2I_s^2 + \tau_{ia,r}^2I_a^2 + 2\tau_{is,r}\tau_{ia,r}I_sI_a + \tau_{r,s}^2R^2 \\~\\
	
	&& + 2(\tau_{is,r}I_s + \tau_{ia,r}I_a)(-\tau_{r,s}R) + \tau_{h,d}^2H^2 + \tau_{is,d}^2I_s^2\\~\\
	
	&& + 2\tau_{h,d}\tau_{is,d}HI_s\\~\\

	&\leq& 2(\zeta_{ia,s}I_a + \zeta_{is,s}I_s + \zeta_{h,s}H)^2 + \tau_{r,s}^2R^2 + 2\tau_{e,ia}^2E^2 + (\tau_{ia,is}+\tau_{ia,r})^2I_a^2 + \tau_{ia,is}^2I_a^2 \\~\\
	
	&& + (\tau_{is,r}+\tau_{is,h}+\tau_{is,d})^2I_s^2 + \tau_{is,h}^2I_s^2 + \tau_{h,d}^2H^2 + \tau_{is,r}^2I_s^2 + \tau_{ia,r}^2I_a^2 \\~\\
	
	&& + 2\tau_{is,r}\tau_{ia,r}I_sI_a + \tau_{r,s}^2R^2 + 2(\tau_{is,r}I_s + \tau_{ia,r}I_a)(-\tau_{r,s}R) \\~\\
	
	&& + \tau_{h,d}^2H^2 + \tau_{is,d}^2I_s^2 + 2\tau_{h,d}\tau_{is,d}HI_s
	
\end{array}$$

$$\begin{array}{rlll}
	
	&=& 2(\zeta_{ia,s}^2I_a^2 + \zeta_{is,s}^2I_s^2 + 2\zeta_{ia,s}\zeta_{is,s}I_aI_s + \zeta_{h,s}^2H^2 + 2\zeta_{h,s}\zeta_{ia,s}HI_a + 2\zeta_{h,s}\zeta_{is,s}HI_s)\\~\\ 
	
	&& + \tau_{r,s}^2R^2 + 2\tau_{e,ia}^2E^2 + (\tau_{ia,is}+\tau_{ia,r})^2I_a^2 + \tau_{ia,is}^2I_a^2 + (\tau_{is,r}+\tau_{is,h}+\tau_{is,d})^2I_s^2 \\~\\
	
	&& + \tau_{is,h}^2I_s^2 + \tau_{h,d}^2H^2 + \tau_{is,r}^2I_s^2 + \tau_{ia,r}^2I_a^2 + 2\tau_{is,r}\tau_{ia,r}I_sI_a + \tau_{r,s}^2R^2 \\~\\
	
	&& + 2(-\tau_{is,r}\tau_{r,s}I_sR - \tau_{ia,r}\tau_{r,s}I_aR) + \tau_{h,d}^2H^2 + \tau_{is,d}^2I_s^2 \\~\\
	
	&& + 2\tau_{h,d}\tau_{is,d}HI_s \\~\\
	
	&\leq& 2(\zeta_{ia,s}^2I_a^2 + \zeta_{is,s}^2I_s^2 + 2\zeta_{ia,s}\zeta_{is,s}I_aI_s + \zeta_{h,s}^2H^2 + 2\zeta_{h,s}\zeta_{ia,s}HI_a + 2\zeta_{h,s}\zeta_{is,s}HI_s)\\~\\ 
	&& + \tau_{r,s}^2R^2 + 2\tau_{e,ia}^2E^2 + (\tau_{ia,is}+\tau_{ia,r})^2I_a^2 + \tau_{ia,is}^2I_a^2 + (\tau_{is,r}+\tau_{is,h}+\tau_{is,d})^2I_s^2 \\~\\
	
	&& + \tau_{is,h}^2I_s^2 + \tau_{h,d}^2H^2 + \tau_{is,r}^2I_s^2 + \tau_{ia,r}^2I_a^2 + 2\tau_{is,r}\tau_{ia,r}I_sI_a \\~\\
	
	&& + \tau_{r,s}^2R^2 + \tau_{h,d}^2H^2 + \tau_{is,d}^2I_s^2 + 2\tau_{h,d}\tau_{is,d}HI_s \\~\\
	
	&\leq& \left[ 2(\zeta_{ia,s}^2 + \zeta_{is,s}^2 + 2\zeta_{ia,s}\zeta_{is,s} + \zeta_{h,s}^2 + 2\zeta_{h,s}\zeta_{ia,s} + 2\zeta_{h,s}\zeta_{is,s}) + \tau_{r,s}^2 + 2\tau_{e,ia}^2 + \tau_{ia,is}^2 \right.\\~\\
	
	&& + (\tau_{ia,is} + \tau_{ia,r})^2 + (\tau_{is,r}+\tau_{is,h}+\tau_{is,d})^2 + \tau_{is,h}^2 + \tau_{h,d}^2 + \tau_{is,r}^2 + \tau_{ia,r}^2 \\~\\
	
	&& \left. + 2\tau_{is,r}\tau_{ia,r} + \tau_{r,s}^2 + \tau_{h,d}^2 + \tau_{is,d}^2 + 2\tau_{h,d}\tau_{is,d} \right] N^2 \\~\\
	
	&=& \left[ 2(\zeta_{ia,s} + \zeta_{is,s} + \zeta_{h,s})^2 + 2\tau_{r,s}^2 + 2\tau_{e,ia}^2 + (\tau_{ia,is} + \tau_{ia,r})^2 + \tau_{ia,is}^2 + \tau_{is,h}^2 + \tau_{is,r}^2 \right. \\~\\
	
	&& + (\tau_{is,r}+\tau_{is,h} + \tau_{is,d})^2 + \tau_{h,d}^2 + \tau_{ia,r}^2 + 2\tau_{is,r}\tau_{ia,r} + \tau_{h,d}^2 + \tau_{is,d}^2 \\~\\
	
	&& \left. + 2\tau_{h,d}\tau_{is,d} \right] N^2 \\~\\
	
	&=& \left[ 2(\zeta_{ia,s} + \zeta_{is,s} + \zeta_{h,s})^2 + 2\tau_{r,s}^2 + 2\tau_{e,ia}^2 + (\tau_{ia,is} + \tau_{ia,r})^2 + \tau_{ia,is}^2 + \tau_{is,h}^2 + \tau_{ia,r}^2 \right. \\~\\
	
	&& + (\tau_{is,r}+\tau_{is,h} + \tau_{is,d})^2 + \tau_{h,d}^2 + \tau_{is,r}^2 + (\tau_{h,d} + \tau_{is,d})^2 \\~\\
	
	&& \left. + 2\tau_{ia,r}\tau_{is,r} \right] N^2 \\~\\
	
	&=& \left[(\tau_{ia,is} + \tau_{ia,r})^2  + \tau_{h,d}^2 + (\tau_{h,d} + \tau_{is,d})^2 + (\tau_{is,r}+\tau_{is,h} + \tau_{is,d})^2  \right. \\~\\
	
	&& + 2(\zeta_{ia,s} + \zeta_{is,s} + \zeta_{h,s})^2 + \tau_{is,r} + (\tau_{is,r} + 2\tau_{ia,r}) + \tau_{ia,is}^2 + \tau_{is,h}^2 \\~\\
	
	&& \left. + \tau_{ia,r}^2 + 2\tau_{r,s}^2 + 2\tau_{e,ia}^2 \right] N^2
\end{array}$$
This implies that
$$\begin{array}{rlll}
     \|f_1(t,x)\|
	&\leq& N \left[(\tau_{ia,is} + \tau_{ia,r})^2  +\tau_{h,d}^2 + (\tau_{h,d} + \tau_{is,d})^2 + (\tau_{is,r}+\tau_{is,h} + \tau_{is,d})^2  \right. \\~\\
	
	&& + 2(\zeta_{ia,s} + \zeta_{is,s} + \zeta_{h,s})^2 + \tau_{is,r}(\tau_{is,r} + 2\tau_{ia,r}) + \tau_{ia,is}^2 + \tau_{is,h}^2 \\~\\
	
	&& \left. + \tau_{ia,r}^2 + 2\tau_{r,s}^2 + 2\tau_{e,ia}^2 \right]^{1/2}.
\end{array}$$

\noindent
We have also:
$$\begin{array}{rlll}
	&& \|f_2(t,x)\|^2 \\~\\
	
	&\leq& 7(\zeta_{ia,s}I_a+\zeta_{is,s}I_s+\zeta_{h,s}H)^2 + 3(\zeta_{ia,s}I_a+\zeta_{is,s}I_s+\zeta_{h,s}H-\tau_{r,s}R)^2 + 2\tau_{e,ia}^2E^2 \\~\\
	
	&& + (\tau_{ia,is}+\tau_{ia,r})^2I_a^2 -2\tau_{e,ia}(\tau_{ia,is}+\tau_{ia,r})EI_a + 2(\tau_{is,r}+\tau_{is,h}+\tau_{is,d})^2I_s^2 \\~\\
	
	&& + \tau_{ia,is}^2I_a^2  + 2\tau_{is,h}^2I_s^2 + (\tau_{h,d}-\tau_{h,r})^2H^2 + 3\tau_{r,s}^2R^2 + \tau_{ia,r}^2I_a^2 + 2\tau_{is,r}^2I_s^2 + \tau_{h,r}^2H^2 \\~\\
	
	&& + 2\tau_{is,d}^2I_s^2 + \tau_{h,d}^2H^2 \\~\\
	
	&\leq& \left[ 7(\zeta_{ia,s}+\zeta_{is,s}+\zeta_{h,s})^2 + 3(\zeta_{ia,s}+\zeta_{is,s}+\zeta_{h,s})^2 + (\tau_{ia,is}+\tau_{ia,r})^2 + 2\tau_{e,ia}^2 \right. \\~\\
	
	&& + 2(\tau_{is,r}+\tau_{is,h}+\tau_{is,d})^2 + \tau_{ia,is}^2  + 2\tau_{is,h}^2 + (\tau_{h,d}-\tau_{h,r})^2 + 3\tau_{r,s}^2 + \tau_{ia,r}^2 + 2\tau_{is,r}^2 \\~\\
	
	&& \left. + \tau_{h,r}^2 + 2\tau_{is,d}^2 + \tau_{h,d}^2 \right] N^2
\end{array}$$

$$\begin{array}{rlll}
	&=& \left[ 10(\zeta_{ia,s}+\zeta_{is,s}+\zeta_{h,s})^2 + (\tau_{ia,is}+\tau_{ia,r})^2 + 2\tau_{e,ia}^2 \right. \\~\\
	
	&& + 2(\tau_{is,r}+\tau_{is,h}+\tau_{is,d})^2 + \tau_{ia,is}^2  + 2\tau_{is,h}^2 + (\tau_{h,d}-\tau_{h,r})^2 + 3\tau_{r,s}^2 + \tau_{ia,r}^2 + 2\tau_{is,r}^2 \\~\\
	
	&& \left. + \tau_{h,r}^2 + 2\tau_{is,d}^2 + \tau_{h,d}^2 \right] N^2
\end{array}$$

\noindent
This implies that

$$\begin{array}{rlll}
	&& \|f_2(t,x)\| \\~\\
	
	&=& N \left[ 10(\zeta_{ia,s}+\zeta_{is,s}+\zeta_{h,s})^2 + (\tau_{ia,is}+\tau_{ia,r})^2 + 2\tau_{e,ia}^2 \right. \\~\\
	
	&& + 2(\tau_{is,r}+\tau_{is,h}+\tau_{is,d})^2 + \tau_{ia,is}^2  + 2\tau_{is,h}^2 + (\tau_{h,d}-\tau_{h,r})^2 + 3\tau_{r,s}^2 + \tau_{ia,r}^2 + 2\tau_{is,r}^2 \\~\\
	
	&& \left. + \tau_{h,r}^2 + 2\tau_{is,d}^2 + \tau_{h,d}^2 \right]^{1/2}.
\end{array}$$
Hence, $f_0(t,x,u)$ is bounded.

\begin{itemize}
	\item[(iii)] \textbf{Prove that the integrand of $J_1$ is convex with respect to the control}\\
	The integrant of $J_1$ is the Lagrangian defined by
	$$\L_1(t,x,u) = \lambda_1 E + \lambda_2 I_a + \lambda_3 I_s + \lambda_4 H + \dis\frac{1}{2}\dis\sum_{i=1}^{9}b_iu_i,$$
	where $x= (S, E, I_a, I_s, H, R, D)$ represents the \textbf{state variable}.
	
	Let $v=(v_1, v_2, \cdots, v_9), w=(w_1, w_2, \cdots, w_9) \in \U$ and $\alpha \in [0, u_{max}]$. We have:
	$$
	\begin{array}{rllll}
		\L_1(t,x, \alpha v + (1-\alpha)w) &=& \lambda_1 E + \lambda_2 I_a + \lambda_3 I_s + \lambda_4 H + \dis\frac{1}{2}\dis\sum_{i=1}^{9}b_i(\alpha v_i + (1-\alpha)w_i)^2, \\
		\alpha \L_1(t,x,v) &=& \alpha \left[ \lambda_1 E + \lambda_2 I_a + \lambda_3 I_s + \lambda_4 H + \dis\frac{1}{2}\dis\sum_{i=1}^{9}b_iv_i^2 \right], \\
		(1-\alpha) \L_1(t,x,w) &=& (1-\alpha) \left[ \lambda_1 E + \lambda_2 I_a + \lambda_3 I_s + \lambda_4 H + \dis\frac{1}{2}\dis\sum_{i=1}^{9}b_iw_i^2 \right].
	\end{array}
	$$
\end{itemize}

\noindent
This implies that
\begin{equation}\label{conv1}
	\alpha \L_1(t,x,v) + (1-\alpha) \L_1(t,x,w) = \lambda_1 E + \lambda_2 I_a + \lambda_3 I_s + \lambda_4 H + \dis\frac{1}{2}\dis\sum_{i=1}^{9}b_i(\alpha v_i^2 + (1-\alpha)w_i^2)
\end{equation}
and
$$\begin{array}{rlllll}
	&& \L_1(t,x, \alpha v + (1-\alpha)w) - \left\{ \alpha \L_1(t,x,v) + (1-\alpha) \L_1(t,x,w) \right\} \\
	&=& \dis\frac{1}{2} \dis\sum_{i=1}^{9}b_i \left\{\left[\alpha v_i + (1-\alpha)w_i\right]^2 - \left[\alpha v_i^2 + (1-\alpha)w_i^2\right]\right\} \\
	&=& \dis\frac{1}{2} \dis\sum_{i=1}^{9}b_i \left\{ \alpha^2 v_i^2 + (1-\alpha)^2 w_i^2 + 2\alpha(1-\alpha)v_iw_i - \alpha v_i^2 - (1-\alpha)w_i^2 \right\} \\
	&=& \dis\frac{1}{2} \dis\sum_{i=1}^{9}b_i \left\{\alpha(\alpha-1)v_i^2 + (1-\alpha)(1-\alpha-1)w_i^2 + 2\alpha(1-\alpha)v_iw_i \right\} \\
	&=& \dis\frac{1}{2} \dis\sum_{i=1}^{9}b_i \left\{\alpha(\alpha-1)v_i^2 + \alpha(\alpha-1)w_i^2 - 2\alpha(\alpha-1)v_iw_i \right\}
\end{array}$$

So,
$$\begin{array}{rlllll}
	&& \L_1(t,x, \alpha v + (1-\alpha)w) - \left\{ \alpha \L_1(t,x,v) + (1-\alpha) \L_1(t,x,w) \right\} \\
	&=& \dis\frac{1}{2}\alpha(\alpha-1)\dis\sum_{i=1}^{9}b_i \left(v_i^2 + w_i^2 - 2v_iw_i\right) \\
	&=& \dis\frac{1}{2}\alpha (\alpha-1)\dis\sum_{i=1}^{9}b_i (v_i-w_i)^2 \\
	&\leq& 0.
\end{array}$$
Hence, the Lagrangian $\L_1(t,x,u)$ is \textbf{convex}.

\begin{itemize}
	\item[(iv)] \textbf{Prove that the Lagrangian $\L_1$ is bounded below.}\\
	It is clear that
	$$\L_1(t,x,u) = \lambda_1E + \lambda_2I_a + \lambda_3 I_s + \lambda_4 H + \dis\frac{1}{2}\dis\sum_{i=1}^{9}b_iu_i^2 \geq \dis\frac{1}{2}\dis\sum_{i=1}^{9}b_iu_i^2.$$
	
	In addition, each 
	$b_i \geq \dis\min_{i=1,\cdots,9}\left\{b_i\right\}$. So, \\
	$$\begin{array}{rlllll}
		\dis\frac{1}{2}\dis\sum_{i=1}^{9}b_iu_i^2 &\geq& \dis\frac{1}{2}\dis\min_{i=1,\cdots,9}\left\{b_i\right\}\dis\sum_{i=1}^{9}u_i^2 \\
		&\geq& \dis\frac{1}{2}\dis\min_{i=1,\cdots,9}\left\{b_i\right\}\left(\dis\sum_{i=1}^{9}u_i^2\right)^{2/2} - a_1, ~~ a_1 > 0 \\
		&=& a_0 \left(\dis\sum_{i=1}^{9}u_i^2\right)^{a_2/2} - a_1, ~~ \forall a_1 > 0 
	\end{array}$$
	where $a_0 = \dis\frac{1}{2}\dis\min_{i=1,\cdots,9}\left\{b_i\right\}$ and $a_2 = 2$. Hence, according to Proposition \ref{ExistContr}, there exists an optimal control $u^*$ satisfying \eqref{optimiz1}.
\end{itemize}
\end{proof}

\noindent
Now, we characterize the optimal control, by means of  Pontryagin's Maximum Principle, which converts the non-autonomous control system \eqref{contsyst1} (with the objective functional \eqref{objective1}), into a pointwise minimization problem, with respect to the control $u\in \U$; the objective here being to minimize the cost function $J_1$.\\

\noindent
Let $\alpha_i, ~ i=1, \cdots,7$ be the adjoint variables corresponding to the state variables $S, E, I_a, I_s, H, R$ and $D$ respectively. The Hamiltonian $\H$ associated with our problem is defined by: 

\begin{equation}\label{hamilton1}
	\begin{array}{rlllll}
		\H &=& \alpha_1\dis\frac{dS}{dt} + \alpha_2\dis\frac{dE}{dt} + \alpha_3\dis\frac{dI_a}{dt} + \alpha_4\dis\frac{dI_s}{dt} + \alpha_5\dis\frac{dH}{dt} + \alpha_6\dis\frac{dR}{dt} + \alpha_7\dis\frac{dD}{dt} \\
		&& + \left(\lambda_1E + \lambda_2I_a + \lambda_3 I_s + \lambda_4 H + \dis\frac{1}{2}\dis\sum_{i=1}^{9}b_iu_i^2\right),
	\end{array}
\end{equation}
and we have the following result.

\begin{theo}\label{characTheoContr1}
	The optimal control $u^*=(u_1^*,u_2^*, \cdots, u_9^*)$ solution of the optimization problem
	\begin{equation*}
		J_1(u^*) = \min_{u\in\U}{J_1(u)}
	\end{equation*}
	is characterized by:
	\begin{equation}\label{characSyst}
		u_i^* = \left\{
		\begin{array}{clllll}
			0 &\text{if}& \omega_i^* \leq 0 \\
			\omega_i^* &\text{if}& 0<\omega_i<u_{i_{max}}\\
			u_{i_{max}} &\text{if}& u_{i_{max}} \leq \omega_i
		\end{array}
		\right., ~~ i=1, \cdots, 9
	\end{equation}
	where:
	$$\begin{array}{rllllll}
		\omega_1^* &=& \dis\frac{S^*(\alpha_2-\alpha_1)(\zeta_{ia,s}I_a^* + \zeta_{is,s}I_s^* + \zeta_{h,s}H^*)}{b_1N} \\~\\
		
		\omega_2^* &=& \dis\frac{S^*(\alpha_1-\alpha_2)(\zeta_{ia,s}I_a^* + \zeta_{is,s}I_s^* + \zeta_{h,s}H^*)}{b_2N}
	\end{array}$$
	
	$$\begin{array}{rllllll}
		\omega_3^* &=& \dis\frac{S^*(\alpha_2-\alpha_1)(\zeta_{ia,s}I_a^* + \zeta_{is,s}I_s^* + \zeta_{h,s}H^*)}{b_3N} + \dis\frac{(\alpha_1-\alpha_6)\tau_{r,s}R^*}{b_3}\\~\\
		
		\omega_4^* &=& \dis\frac{S^*(\alpha_2-\alpha_1)(\zeta_{ia,s}I_a^* + \zeta_{is,s}I_s^* + \zeta_{h,s}H^*)}{b_4N} + \dis\frac{(\alpha_1-\alpha_6)\tau_{r,s}R^*}{b_4}\\~\\
		
		\omega_5^* &=& \dis\frac{S^*(\alpha_2-\alpha_1)(\zeta_{ia,s}I_a^* + \zeta_{is,s}I_s^* + \zeta_{h,s}H^*)}{b_5N} + \dis\frac{(\alpha_1-\alpha_6)\tau_{r,s}R^*}{b_5}\\~\\
		
		\omega_6^* &=& \dis\frac{(\alpha_3-\alpha_2)\tau_{e,ia}E^* + (\alpha_4-\alpha_3)\tau_{ia,is}I_a^* + (\alpha_6-\alpha_3)\tau_{ia,r}I_a^*}{b_6}\\~\\
		
		\omega_7^* &=& \dis\frac{\left[(\alpha_5\tau_{is,h}+\alpha_6\tau_{is,r}+\alpha_7\tau_{is,d}) - \alpha_4(\tau_{is,r}+\tau_{is,h}+\tau_{is,d})\right]I_s^*}{b_7}
	\end{array}$$
	
	$$\begin{array}{rllllll}
		\omega_8^* &=& \dis\frac{\left[(\alpha_5\tau_{is,h}+\alpha_6\tau_{is,r}+\alpha_7\tau_{is,d}) - \alpha_4(\tau_{is,r}+\tau_{is,h}+\tau_{is,d})\right]I_s^*}{b_8}\\~\\
		
		\omega_9^* &=& \dis\frac{\left[(\alpha_5-\alpha_6)\tau_{h,r} + (\alpha_7-\alpha_5)\tau_{h,d}\right]H^*}{b_9}
	\end{array}$$
	with the septuple $(S^*, E^*, I_a^*, I_s^*, H^*, R^*, D^*)$ solution of the control system \eqref{contsyst1} and the adjoint variables $(\alpha_1, \alpha_2, \cdots, \alpha_7)$ satisfy the following adjoint system, with transversality conditions $\alpha_i(T) = 0, ~~ i=1,\cdots,7$:
	\begin{equation}\label{adjointSyst1}
		\left\{
		\begin{array}{lllllll}
			\dis\frac{d\alpha_1}{dt} &=& \dis\frac{(\alpha_1-\alpha_2)}{N}(1-u_1+u_2-u_3-u_4-u_5)(\zeta_{ia,s}I_a + \zeta_{is,s}I_s + \zeta_{h,s}H) \\~\\		
			\dis\frac{d\alpha_2}{dt} &=& (\alpha_2-\alpha_3)(1-u_6)\tau_{e,ia} -\lambda_1\\~\\
			\dis\frac{d\alpha_3}{dt} &=& \dis\frac{(\alpha_1-\alpha_2)S}{N}(1-u_1+u_2-u_3-u_4-u_5)\zeta_{ia,s} + (1-u_6)\left[(\alpha_3-\alpha_4)\tau_{ia,is} \right.\\ 
			&& \left. + (\alpha_3-\alpha_6)\tau_{ia,r}\right] - \lambda_2 \\~\\
			\dis\frac{d\alpha_4}{dt} &=& \dis\frac{(\alpha_1-\alpha_2)S}{N}(1-u_1+u_2-u_3-u_4-u_5)\zeta_{is,s} + (1-u_7-u_8)\left[(\alpha_4-\alpha_6)\tau_{is,r} \right.\\ 
			&& \left. + (\alpha_4-\alpha_5)\tau_{is,h} + (\alpha_4-\alpha_7)\tau_{is,d}\right] - \lambda_3\\~\\
			\dis\frac{d\alpha_5}{dt} &=& \dis\frac{(\alpha_1-\alpha_2)S}{N}(1-u_1+u_2-u_3-u_4-u_5)\zeta_{h,s} + (1-u_9)(\alpha_5-\alpha_7)\tau_{h,d} \\ 
			&& + u_9(\alpha_5-\alpha_6)\tau_{h,r} - \lambda_4\\~\\
			\dis\frac{d\alpha_6}{dt} &=& \dis(\alpha_6-\alpha_1)(1-u_3-u_4-u_5)\tau_{r,s}\\~\\
			\dis\frac{d\alpha_7}{dt} &=& 0.
		\end{array}
		\right.
	\end{equation}
\end{theo}

\noindent
\begin{proof}
Consider the following Hamilton's equations, derived from the total energy of the system (called Hamiltonian):\\
$\dis\frac{\partial \H}{\partial S} = -\dis\frac{d\alpha_1}{dt},
	\dis\frac{\partial \H}{\partial E} = -\dis\frac{d\alpha_2}{dt},
	\dis\frac{\partial \H}{\partial I_a} = -\dis\frac{d\alpha_3}{dt},
	\dis\frac{\partial \H}{\partial I_s} = -\dis\frac{d\alpha_4}{dt},
	\dis\frac{\partial \H}{\partial H} = -\dis\frac{d\alpha_5}{dt},
	\dis\frac{\partial \H}{\partial R} = -\dis\frac{d\alpha_6}{dt},
	\dis\frac{\partial \H}{\partial D} = -\dis\frac{d\alpha_7}{dt}
$.\\

\noindent
The expanded form of the Hamiltonian is given below:
\begin{eqnarray*}
	\H &=& \alpha_1 \left[ -\left(1-u_1+u_2-u_3-u_4-u_5\right) \left(\zeta_{ia,s}I_a + \zeta_{is,s}I_s + \zeta_{h,s}H\right)\dis\frac{S}{N} + (1-u_3-u_4-u_5)\tau_{r,s}R \right] \\
	&& + \alpha_2 \left[ \left(1-u_1+u_2-u_3-u_4-u_5\right) \left(\zeta_{ia,s}I_a + \zeta_{is,s}I_s + \zeta_{h,s}H\right)\dis\frac{S}{N} - (1-u_6)\tau_{e,ia}E \right] \\
	&& + \alpha_3 \left[ (1-u_6)\tau_{e,ia}E - (1-u_6) (\tau_{ia,is} + \tau_{ia,r})I_a \right]\\
	&& + \alpha_4 \left[ (1-u_6)\tau_{ia,is}I_a - (1-u_7-u_8)(\tau_{is,r} + \tau_{is,h} + \tau_{is,d})I_s \right] \\
	&& + \alpha_5 \left[ (1-u_7-u_8)\tau_{is,h}I_s - (1-u_9)\tau_{h,d}H -u_9\tau_{h,r}H\right] \\
	&& + \alpha_6 \left[ (1-u_7-u_8)\tau_{is,r}I_s + (1-u_6)\tau_{ia,r}I_a + u_9\tau_{h,r}H \right. \\
	&& \left. -(1-u_3-u_4-u_5)\tau_{r,s}R \right] + \alpha_7 \left[ (1-u_9)\tau_{h,d}H + (1-u_7-u_8)\tau_{is,d}I_s \right]\\
	&& + \lambda_1E + \lambda_2I_a + \lambda_3 I_s + \lambda_4 H + \dis\frac{1}{2}\dis\sum_{i=1}^{9}b_iu_i^2.
\end{eqnarray*}

\noindent
To construct the adjoint state, let's differentiate the Hamiltonian $\H$ with respect to each state variable. We have:
$$\dis\frac{\partial \H}{\partial S} = (\alpha_2-\alpha_1)(1-u_1+u_2-u_3-u_4-u_5)(\zeta_{ia,s}I_a + \zeta_{is,s}I_s + \zeta_{h,s}H)\dis\frac{1}{N}$$
which implies that
$$\dis\frac{d \alpha_1}{dt} = \dis\frac{(\alpha_1-\alpha_2)}{N}(1-u_1+u_2-u_3-u_4-u_5)(\zeta_{ia,s}I_a + \zeta_{is,s}I_s + \zeta_{h,s}H).$$

\noindent
Similarly, we have the following remaining equations of the adjoint system:

\begin{eqnarray*}
	\dis\frac{d \alpha_2}{dt} &=& (\alpha_2-\alpha_3)(1-u_6)\tau_{e,ia} - \lambda_1 \\
	\dis\frac{d \alpha_3}{dt} &=& \dis\frac{(\alpha_1-\alpha_2)S}{N}(1-u_1+u_2-u_3-u_4-u_5)\zeta_{ia,s} + (1-u_6)\left[(\alpha_3-\alpha_4)\tau_{ia,is}\right.\\
	&& \left. + (\alpha_3-\alpha_6)\tau_{ia,r}\right] - \lambda_2\\
	\dis\frac{d \alpha_4}{dt} &=& \dis\frac{(\alpha_1-\alpha_2)S}{N}(1-u_1+u_2-u_3-u_4-u_5)\zeta_{is,s} + (1-u_7-u_8)\left[(\alpha_4-\alpha_6)\tau_{is,r}\right.\\
	&& \left. + (\alpha_4-\alpha_5)\tau_{is,h} + (\alpha_4-\alpha_7)\tau_{is,d}\right] - \lambda_3\\
	\dis\frac{d \alpha_5}{dt} &=& \dis\frac{(\alpha_1-\alpha_2)S}{N}(1-u_1+u_2-u_3-u_4-u_5)\zeta_{h,s} + (1-u_9)(\alpha_5-\alpha_7)\tau_{h,d}\\
	&& + u_9(\alpha_5-\alpha_6)\tau_{h,r} - \lambda_4 \\
	\dis\frac{d \alpha_6}{dt} &=& (\alpha_6-\alpha_1)(1-u_3-u_4-u_5)\tau_{r,s}\\
	\dis\frac{d \alpha_7}{dt} &=& 0.
\end{eqnarray*}

\begin{rem}\label{Tranversality}
	Since the time horizon $T$ is fixed and there is no explicit dependency on the state variable at $T$, then our optimization problem does not require any feedback from the terminal state. Consequently, the adjoint variable at $T$ will satisfy: $\alpha_i(T) = 0, ~~ i=1,\cdots, 7$.
\end{rem}

\noindent
Now, for each $u_i,~~ i=1,\cdots,9$, we need to solve the following
$$\dis\frac{\partial \H}{\partial u_i} = 0$$
to complete the characterization of the optimal control.

\noindent
We have:\\
$$\dis\frac{\partial \H}{\partial u_1} = (\alpha_1-\alpha_2)(\zeta_{ia,s}I_a + \zeta_{is,s}I_s + \zeta_{h,s}H)\dis\frac{S}{N} + b_1u_1.$$

\noindent
Solving $\dis\frac{\partial \H}{\partial u_1} = 0$ for $u_1$, we obtain:

$$u_1 = \dis\frac{S(\alpha_2-\alpha_1)(\zeta_{ia,s}I_a + \zeta_{is,s}I_s + \zeta_{h,s}H)}{b_1N}.$$

\noindent
Taking into account the admissible control set, if we set $u_{1_{max}} = \dis\max\{u_1\}$ for all $u_1\in [0,u_{max}]$, then it is clear that the optimal value of the control variable $u_1$ is given by:
$$u_1^* = \min \left\{\max\left\{0, \omega_1^*\right\}, u_{1_{max}}\right\},$$
where
$$\omega_1^* = \dis\frac{S^*(\alpha_2-\alpha_1)(\zeta_{ia,s}I_a^* + \zeta_{is,s}I_s^* + \zeta_{h,s}H^*)}{b_1N}$$
Hence:

$$u_1^* = \left\{
\begin{array}{rlll}
	0 &\text{if}& \omega_1^* \leq 0 \\
	\omega_1^* &\text{if}& 0 < \omega_1^* < u_{1_{max}} \\
	u_{1_{max}} &\text{if}& u_{1_{max}} \leq \omega_1^*
\end{array}
\right.$$
Similarly, we have:

$$u_i^* = \left\{
\begin{array}{rlll}
	0 &\text{if}& \omega_i^* \leq 0 \\
	\omega_i^* &\text{if}& 0 < \omega_i^* < u_{i_{max}} \\
	u_{i_{max}} &\text{if}& u_{i_{max}} \leq \omega_i^*
\end{array}
\right., ~~ i=2,\cdots, 9$$
where:

\begin{eqnarray*}
	\omega_2^* &=& \dis\frac{S^*(\alpha_1-\alpha_2)(\zeta_{ia,s}I_a^* + \zeta_{is,s}I_s^* + \zeta_{h,s}H^*)}{b_2N}\\~\\
	\omega_3^* &=& \dis\frac{S^*(\alpha_2-\alpha_1)(\zeta_{ia,s}I_a^* + \zeta_{is,s}I_s^* + \zeta_{h,s}H^*)}{b_3N} + \dis\frac{(\alpha_1-\alpha_6)\tau_{r,s}R^*}{b_3}\\~\\
	\omega_4^* &=& \dis\frac{S^*(\alpha_2-\alpha_1)(\zeta_{ia,s}I_a^* + \zeta_{is,s}I_s^* + \zeta_{h,s}H^*)}{b_4N} + \dis\frac{(\alpha_1-\alpha_6)\tau_{r,s}R^*}{b_4}\\~\\
	\omega_5^* &=& \dis\frac{S^*(\alpha_2-\alpha_1)(\zeta_{ia,s}I_a^* + \zeta_{is,s}I_s^* + \zeta_{h,s}H^*)}{b_5N} + \dis\frac{(\alpha_1-\alpha_6)\tau_{r,s}R^*}{b_5}\\~\\
	\omega_6^* &=& \dis\frac{(\alpha_3-\alpha_2)\tau_{e,ia}E^* + (\alpha_4-\alpha_3)\tau_{ia,is}I_a^* + (\alpha_6-\alpha_3)\tau_{ia,r}I_a^*}{b_6}
\end{eqnarray*}

\begin{eqnarray*}
	\omega_7^* &=& \dis\frac{\left[(\alpha_5\tau_{is,h} + \alpha_6\tau_{is,r} + \alpha_7\tau_{is,d}) - \alpha_4(\tau_{is,r} + \tau_{is,h} + \tau_{is,d})\right]I_s^*}{b_7} \\~\\
	\omega_8^* &=& \dis\frac{\left[(\alpha_5\tau_{is,h} + \alpha_6\tau_{is,r} + \alpha_7\tau_{is,d}) - \alpha_4(\tau_{is,r} + \tau_{is,h} + \tau_{is,d})\right]I_s^*}{b_8} \\~\\
	\omega_9^* &=& \dis\frac{\left[(\alpha_5-\alpha_6)\tau_{h,r} + (\alpha_7 - \alpha_5)\tau_{h,d}\right]H^*}{b_9}.
\end{eqnarray*}

\noindent
This ends the proof of Theorem \ref{characTheoContr1}.
\end{proof}

\subsection{Effectiveness Objective Function}
In this section, we aim to assess the success of intervention strategies in reducing infections, decreasing hospitalizations and mortality, and ultimately controlling the epidemic. Effective interventions would minimize transmission, protect vulnerable populations, and support a steady decline in outbreak intensity. In other words, the objective is to maximize effectiveness, as achieving high efficacy is crucial to containing the virus's spread. Let's consider the effectiveness objective function
\begin{equation}\label{objective2}
	\small
	J_2(u_1, u_2, \cdots, u_9) = \dis\int_{0}^{T} \left[ (b_1u_1 + b_2u_2)S + \dis\sum_{i=3}^{5}b_iu_i(S+R) + b_6u_6(E+I_a) + (b_7u_7 + b_8u_8)I_s + b_9u_9H\right]e^{-\sigma t}dt,
\end{equation}
where $\left\{b_i\right\}_{i=1}^{9}$ are the balancing cost factors associated with the different controls $\left\{u_i\right\}_{i=1}^{9}$; $\sigma$ is a scaling parameter and $T$ is the final time for control implementation.

\noindent
Our goal is to find an optimal control $u^* = (u_1^*, u_2^*, \cdots, u_9^*)$ such that 
\begin{equation}
	J_2(u^*) = \max_{u\in\U}{J_2(u)}
\end{equation}
where $u=(u_1,u_2, \cdots, u_9) \in \U$.

\begin{lem}\label{convex2}
	The integrand of the objective function \eqref{objective2} is convex with respect to the control $u$.
\end{lem}

\noindent
\begin{proof}
	The integrand of \eqref{objective2} is the Lagrangian denoted by
	$$\mathcal{L}_2(t,x,u) = \left((b_1u_1 + b_2u_2)S + \dis\sum_{i=3}^{5}b_iu_i(S+R) + b_6u_6(E+I_a) + (b_7u_7 + b_8u_8)I_s + b_9u_9H\right)e^{-\sigma t}.$$
	
	\noindent
	Let $v,w\in \U, ~~ \alpha\in [0, u_{max}]$. We have:
	
	$$\begin{array}{rlllll}
		&& \L(t,x,\alpha v + (1-\alpha)w) \\~\\
		&=& \left\{ [ b_1 (\alpha v_1 + (1-\alpha)w_1) + b_2 (\alpha v_2 + (1-\alpha)w_2)]S + \dis\sum_{i=3}^{5}b_i(\alpha v_i + (1-\alpha)w_i)(S+R) \right.\\~\\
		&& + b_6(\alpha v_6 + (1-\alpha)w_6)(E+I_a) + [b_7(\alpha v_7 + (1-\alpha)w_7) + b_8 (\alpha v_8 + (1-\alpha)w_8)]I_s \\~\\
		&& \left. b_9(\alpha v_9 + (1-\alpha)w_9)H \right\}e^{-\sigma t}
	\end{array}$$
	
	$$\begin{array}{rlll}
		&=& \alpha \left\{ (b_1 v_1 +b_2v_2)S + \dis\sum_{i=3}^{5}b_iv_i(S+R) + b_6v_6(E+I_a) + (b_7v_7 + b_8v_8)I_s \right\}e^{-\sigma t} \\~\\
		&& + (1-\alpha)\left\{(b_1w_1 + b_2w_2)S + \dis\sum_{i=3}^{5}b_iw_i(S+R) + b_6w_6(E+I_a) + (b_7w_7 + b_8w_8)I_s\right\}e^{-\sigma t} \\~\\
		&=& \alpha\L (t,x,v) + (1-\alpha)\L(t,x,w).
	\end{array}$$
	Hence, $\L$ is convex.\\
\end{proof}

\noindent
Considering Lemma \ref{convex2} and Theorem \eqref{control1}, it is clear that there exists an optimal control $u^* = (u_1^*, u_2^*, \cdots, u_9^*)$ such that 
\begin{equation}\label{optimiz2}
	J_2(u^*) = \max_{u\in\U}{J_2(u)}.
\end{equation}

\noindent
Now, by means of  Pontryagin's Maximum Principle, we characterize the optimal control $u^*$, solution of \eqref{optimiz2}. The Hamiltonian $\H$ associated with our problem is defined by: 

\begin{equation}\label{hamilton1}
	\begin{array}{rlllll}
		\H &=& \alpha_1\dis\frac{dS}{dt} + \alpha_2\dis\frac{dE}{dt} + \alpha_3\dis\frac{dI_a}{dt} + \alpha_4\dis\frac{dI_s}{dt} + \alpha_5\dis\frac{dH}{dt} + \alpha_6\dis\frac{dR}{dt} + \alpha_7\dis\frac{dD}{dt} \\
		&& + \left((b_1u_1 + b_2u_2)S + \dis\sum_{i=3}^{5}b_iu_i(S+R) + b_6u_6(E+I_a) + (b_7u_7 + b_8u_8)I_s + b_9u_9H\right)e^{-\sigma t},
	\end{array}
\end{equation}
where $\alpha_i, ~ i=1, \cdots,7$ represent the adjoint variables associated with the state variables $S, E, I_a, I_s, H, R$ and $D$ respectively. We have the following result.

\begin{theo}\label{characTheoContr2}
	The optimal conditions for effectiveness are the following:
	\begin{equation}\label{effectivenessCond} \small
		\left\{
		\begin{array}{llllll}
			\dis\frac{(\alpha_1-\alpha_2)S}{N}(\zeta_{ia,s}I_a + \zeta_{is,s}I_s + \zeta_{h,s}H) = -b_1Se^{-\sigma t} &\text{for}& u_1^*\in (0, u_{max}) \\~\\
			\dis\frac{(\alpha_2-\alpha_1)S}{N}(\zeta_{ia,s}I_a + \zeta_{is,s}I_s + \zeta_{h,s}H) = -b_2Se^{-\sigma t} &\text{for}& u_2^*\in (0, u_{max}) \\~\\
			\dis\frac{(\alpha_1-\alpha_2)S}{N}(\zeta_{ia,s}I_a + \zeta_{is,s}I_s + \zeta_{h,s}H) + (\alpha_6-\alpha_1)\tau_{r,s}R = -b_3(S+R)e^{-\sigma t} &\text{for}& u_3^*\in (0, u_{max}) \\~\\
			\dis\frac{(\alpha_1-\alpha_2)S}{N}(\zeta_{ia,s}I_a + \zeta_{is,s}I_s + \zeta_{h,s}H) + (\alpha_6-\alpha_1)\tau_{r,s}R = -b_4(S+R)e^{-\sigma t} &\text{for}& u_4^*\in (0, u_{max}) \\~\\
			\dis\frac{(\alpha_1-\alpha_2)S}{N}(\zeta_{ia,s}I_a + \zeta_{is,s}I_s + \zeta_{h,s}H) + (\alpha_6-\alpha_1)\tau_{r,s}R = -b_5(S+R)e^{-\sigma t} &\text{for}& u_5^*\in (0, u_{max}) \\~\\
			(\alpha_2-\alpha_3)\tau_{e,ia}E + [(\alpha_3-\alpha_4)\tau_{ia,is}+(\alpha_3-\alpha_6)\tau_{ia,r}]I_a = -b_6(E+I_a)e^{-\sigma t} &\text{for}& u_6^*\in (0, u_{max}) \\~\\
			((\alpha_4-\alpha_5)\tau_{is,h} + (\alpha_4-\alpha_6)\tau_{is,r}+(\alpha_4-\alpha_7)\tau_{is,d})I_s = -b_7I_se^{-\sigma t} &\text{for}& u_7^*\in (0, u_{max})\\~\\
			((\alpha_4-\alpha_5)\tau_{is,h} + (\alpha_4-\alpha_6)\tau_{is,r}+(\alpha_4-\alpha_7)\tau_{is,d})I_s = -b_8I_se^{-\sigma t} &\text{for}& u_8^*\in (0, u_{max})\\~\\
			((\alpha_6-\alpha_5)\tau_{h,r}+(\alpha_5-\alpha_7)\tau_{h,d})H = -b_9He^{-\sigma t} &\text{for}& u_9^*\in (0, u_{max})
		\end{array}
		\right.
	\end{equation}
	
	\noindent
	with the septuple $(S^*, E^*, I_a^*, I_s^*, H^*, R^*, D^*)$ solution of the control system \eqref{contsyst1} and the adjoint variables $(\alpha_1, \alpha_2, \cdots, \alpha_7)$ satisfy the following adjoint system, with transversality conditions $\alpha_i(T) = 0, ~~ i=1,\cdots,7$:
	\begin{equation}\label{adjointSyst2}
		\left\{
		\begin{array}{lllllll}
			\dis\frac{d\alpha_1}{dt} &=& \dis\frac{(\alpha_1-\alpha_2)}{N}(1-u_1+u_2-u_3-u_4-u_5)(\zeta_{ia,s}I_a + \zeta_{is,s}I_s + \zeta_{h,s}H) - \dis\sum_{i=1}^{5}b_iu_ie^{-\sigma t} \\~\\		
			\dis\frac{d\alpha_2}{dt} &=& (\alpha_2-\alpha_3)(1-u_6)\tau_{e,ia} - b_6u_6e^{-\sigma t}\\~\\
			\dis\frac{d\alpha_3}{dt} &=& \dis\frac{(\alpha_1-\alpha_2)S}{N}(1-u_1+u_2-u_3-u_4-u_5)\zeta_{ia,s} + (1-u_6)\left[(\alpha_3-\alpha_4)\tau_{ia,is} \right.\\ 
			&& \left. + (\alpha_3-\alpha_6)\tau_{ia,r}\right] - b_6u_6e^{-\sigma t} \\~\\
			\dis\frac{d\alpha_4}{dt} &=& \dis\frac{(\alpha_1-\alpha_2)S}{N}(1-u_1+u_2-u_3-u_4-u_5)\zeta_{is,s} + (1-u_7-u_8)\left[(\alpha_4-\alpha_6)\tau_{is,r} \right.\\ 
			&& \left. + (\alpha_4-\alpha_5)\tau_{is,h} + (\alpha_4-\alpha_7)\tau_{is,d}\right] - (b_7u_7 + b_8u_8)e^{-\sigma t}\\~\\
			\dis\frac{d\alpha_5}{dt} &=& \dis\frac{(\alpha_1-\alpha_2)S}{N}(1-u_1+u_2-u_3-u_4-u_5)\zeta_{h,s} + (1-u_9)(\alpha_5-\alpha_7)\tau_{h,d} \\ 
			&& + u_9(\alpha_5-\alpha_6)\tau_{h,r} - b_9u_9e^{-\sigma t}\\~\\
			\dis\frac{d\alpha_6}{dt} &=& \dis(\alpha_6-\alpha_1)(1-u_3-u_4-u_5)\tau_{r,s} - \dis\sum_{i=3}^{5}b_iu_ie^{-\sigma t} \\~\\
			\dis\frac{d\alpha_7}{dt} &=& 0.
		\end{array}
		\right.
	\end{equation}
\end{theo}

\noindent
\begin{proof}
	Consider the Hamilton's equations, derived from the total energy of the system, we can find the adjoint system as follows\\
	
	\noindent
	The expanded form of the Hamiltonian is given below:
	
	\begin{eqnarray*}
		\H &=& \alpha_1 \left[ -\left(1-u_1+u_2-u_3-u_4-u_5\right) \left(\zeta_{ia,s}I_a + \zeta_{is,s}I_s + \zeta_{h,s}H\right)\dis\frac{S}{N} + (1-u_3-u_4-u_5)\tau_{r,s}R \right] \\
		&& + \alpha_2 \left[ \left(1-u_1+u_2-u_3-u_4-u_5\right) \left(\zeta_{ia,s}I_a + \zeta_{is,s}I_s + \zeta_{h,s}H\right)\dis\frac{S}{N} - (1-u_6)\tau_{e,ia}E \right] \\
		&& + \alpha_3 \left[ (1-u_6)\tau_{e,ia}E - (1-u_6) (\tau_{ia,is} + \tau_{ia,r})I_a \right]\\
		&& + \alpha_4 \left[ (1-u_6)\tau_{ia,is}I_a - (1-u_7-u_8)(\tau_{is,r} + \tau_{is,h} + \tau_{is,d})I_s \right] \\
		&& + \alpha_5 \left[ (1-u_7-u_8)\tau_{is,h}I_s - (1-u_9)\tau_{h,d}H -u_9\tau_{h,r}H\right] \\
		&& + \alpha_6 \left[ (1-u_7-u_8)\tau_{is,r}I_s + (1-u_6)\tau_{ia,r}I_a + u_9\tau_{h,r}H \right. \\
		&& \left. -(1-u_3-u_4-u_5)\tau_{r,s}R \right] + \alpha_7 \left[ (1-u_9)\tau_{h,d}H + (1-u_7-u_8)\tau_{is,d}I_s \right]\\
		&& + \left((b_1u_1 + b_2u_2)S + \dis\sum_{i=3}^{5}b_iu_i(S+R) + b_6u_6(E+I_a) + (b_7u_7 + b_8u_8)I_s + b_9u_9H\right)e^{-\sigma t}.
	\end{eqnarray*}
	
	\noindent
	To construct the adjoint state, let's differentiate the Hamiltonian $\H$ with respect to each state variable. We have:
	$$\dis\frac{\partial \H}{\partial S} = (\alpha_2-\alpha_1)(1-u_1+u_2-u_3-u_4-u_5)(\zeta_{ia,s}I_a + \zeta_{is,s}I_s + \zeta_{h,s}H)\dis\frac{1}{N} + \dis\sum_{i=1}^{5}b_iu_ie^{-\sigma t}$$
	which implies that
	$$\dis\frac{d \alpha_1}{dt} = \dis\frac{(\alpha_1-\alpha_2)}{N}(1-u_1+u_2-u_3-u_4-u_5)(\zeta_{ia,s}I_a + \zeta_{is,s}I_s + \zeta_{h,s}H) - \dis\sum_{i=1}^{5}b_iu_ie^{-\sigma t}.$$
	
	\noindent
	Similarly, we obtain the remaining equations of the adjoint system \eqref{adjointSyst2}. Moreover, according to Remark \ref{Tranversality} we have $\alpha_i(T) = 0, ~~ i=1,\cdots, 7$.\\
	
	\noindent
	Now, for each $u_i,~~ i=1,\cdots,9$, we need to solve the following
	$$\dis\frac{\partial \H}{\partial u_i} = 0$$
	to complete the characterization of the optimal control.
	
	\noindent
	We have the following:
	\begin{equation}\label{partialHu1}
		\dis\frac{\partial \H}{\partial u_1} = \dis\frac{(\alpha_1-\alpha_2)S}{N}(\zeta_{ia,s}I_a + \zeta_{is,s}I_s + \zeta_{h,s}H) - \left(-b_1Se^{-\sigma t}\right),
	\end{equation}
	where $\dis\frac{(\alpha_1-\alpha_2)S}{N}(\zeta_{ia,s}I_a + \zeta_{is,s}I_s + \zeta_{h,s}H)$ represents the marginal benefit of taking the vaccine and $-b_1Se^{-\sigma t}$ the marginal cost. Similarly, we have:
	
	$$
	\begin{array}{rllll}
		\dis\frac{\partial \H}{\partial u_2} &=& \dis\frac{(\alpha_2-\alpha_1)S}{N}(\zeta_{ia,s}I_a + \zeta_{is,s}I_s + \zeta_{h,s}H) + b_2Se^{-\sigma t}.\\~\\
		
		\dis\frac{\partial \H}{\partial u_3} &=& \dis\frac{(\alpha_1-\alpha_2)S}{N}(\zeta_{ia,s}I_a + \zeta_{is,s}I_s + \zeta_{h,s}H) + (\alpha_6-\alpha_1)\tau_{r,s}R + b_3(S+R)e^{-\sigma t}.\\~\\
		
		\dis\frac{\partial \H}{\partial u_4} &=& \dis\frac{(\alpha_1-\alpha_2)S}{N}(\zeta_{ia,s}I_a + \zeta_{is,s}I_s + \zeta_{h,s}H) + (\alpha_6-\alpha_1)\tau_{r,s}R + b_4(S+R)e^{-\sigma t}.\\~\\
		
		\dis\frac{\partial \H}{\partial u_5} &=& \dis\frac{(\alpha_1-\alpha_2)S}{N}(\zeta_{ia,s}I_a + \zeta_{is,s}I_s + \zeta_{h,s}H) + (\alpha_6-\alpha_1)\tau_{r,s}R + b_5(S+R)e^{-\sigma t}.\\~\\
		
		\dis\frac{\partial \H}{\partial u_6} &=& (\alpha_2-\alpha_3)\tau_{e,ia}E + [(\alpha_3-\alpha_4)\tau_{ia,is} + (\alpha_3-\alpha_6)\tau_{ia,r}]I_a + b_6(E+I_a)e^{-\sigma t}.\\~\\
		
		\dis\frac{\partial \H}{\partial u_7} &=& [(\alpha_4-\alpha_5)\tau_{is,h} + (\alpha_4-\alpha_6)\tau_{is,r} + (\alpha_4-\alpha_7)\tau_{is,d}]I_s + b_7I_se^{-\sigma t}.\\~\\
		
		\dis\frac{\partial \H}{\partial u_8} &=& [(\alpha_4-\alpha_5)\tau_{is,h} + (\alpha_4-\alpha_6)\tau_{is,r} + (\alpha_4-\alpha_7)\tau_{is,d}]I_s + b_8I_se^{-\sigma t}.\\~\\
		
		\dis\frac{\partial \H}{\partial u_9} &=& [(\alpha_6-\alpha_5)\tau_{h,r} + (\alpha_5-\alpha_7)\tau_{h,d}]H + b_9He^{-\sigma t}.
	\end{array}
	$$
	
	\noindent
	Based on equation \eqref{partialHu1} and the fact that $u_1\in [0, u_{max}]$, we have
	
	$$
	u_1^* \left\{
	\begin{array}{lllll}
		= u_{max} &\text{if}& \dis\frac{(\alpha_1-\alpha_2)S}{N}(\zeta_{ia,s}I_a + \zeta_{is,s}I_s + \zeta_{h,s}H) > -b_1Se^{-\sigma t} \\
		\in (0, u_{max}) &\text{if}& \dis\frac{(\alpha_1-\alpha_2)S}{N}(\zeta_{ia,s}I_a + \zeta_{is,s}I_s + \zeta_{h,s}H) = -b_1Se^{-\sigma t}\\
		= 0 &\text{if}& \dis\frac{(\alpha_1-\alpha_2)S}{N}(\zeta_{ia,s}I_a + \zeta_{is,s}I_s + \zeta_{h,s}H) < -b_1Se^{-\sigma t}.
	\end{array}
	\right. 
	$$
	
	\noindent
	Analogously, we obtain the following remaining equations completing the characterization of the control $u^*$:
	
		\begin{eqnarray*}
		u_2^* \left\{
		\begin{array}{lllll}
			= u_{max} &\text{if}& \dis\frac{(\alpha_2-\alpha_1)S}{N}(\zeta_{ia,s}I_a + \zeta_{is,s}I_s + \zeta_{h,s}H) > -b_2Se^{-\sigma t} \\
			\in (0, u_{max}) &\text{if}& \dis\frac{(\alpha_2-\alpha_1)S}{N}(\zeta_{ia,s}I_a + \zeta_{is,s}I_s + \zeta_{h,s}H) = -b_2Se^{-\sigma t}\\
			= 0 &\text{if}& \dis\frac{(\alpha_2-\alpha_1)S}{N}(\zeta_{ia,s}I_a + \zeta_{is,s}I_s + \zeta_{h,s}H) < -b_2Se^{-\sigma t}
		\end{array}
		\right. \\~\\
		u_3^* \left\{
		\begin{array}{lllll}
			= u_{max} &\text{if}& \dis\frac{(\alpha_1-\alpha_2)S}{N}(\zeta_{ia,s}I_a + \zeta_{is,s}I_s + \zeta_{h,s}H) + (\alpha_6-\alpha_1)\tau_{r,s}R > -b_3(S+R)e^{-\sigma t} \\
			\in (0, u_{max}) &\text{if}& \dis\frac{(\alpha_1-\alpha_2)S}{N}(\zeta_{ia,s}I_a + \zeta_{is,s}I_s + \zeta_{h,s}H) + (\alpha_6-\alpha_1)\tau_{r,s}R = -b_3(S+R)e^{-\sigma t}\\
			= 0 &\text{if}& \dis\frac{(\alpha_1-\alpha_2)S}{N}(\zeta_{ia,s}I_a + \zeta_{is,s}I_s + \zeta_{h,s}H) + (\alpha_6-\alpha_1)\tau_{r,s}R < -b_3(S+R)e^{-\sigma t}
		\end{array}
		\right.\\~\\
		u_4^* \left\{
		\begin{array}{lllll}
			= u_{max} &\text{if}& \dis\frac{(\alpha_1-\alpha_2)S}{N}(\zeta_{ia,s}I_a + \zeta_{is,s}I_s + \zeta_{h,s}H) + (\alpha_6-\alpha_1)\tau_{r,s}R > -b_4(S+R)e^{-\sigma t} \\
			\in (0, u_{max}) &\text{if}& \dis\frac{(\alpha_1-\alpha_2)S}{N}(\zeta_{ia,s}I_a + \zeta_{is,s}I_s + \zeta_{h,s}H) + (\alpha_6-\alpha_1)\tau_{r,s}R = -b_4(S+R)e^{-\sigma t}\\
			= 0 &\text{if}& \dis\frac{(\alpha_1-\alpha_2)S}{N}(\zeta_{ia,s}I_a + \zeta_{is,s}I_s + \zeta_{h,s}H) + (\alpha_6-\alpha_1)\tau_{r,s}R < -b_4(S+R)e^{-\sigma t}
		\end{array}
		\right.\\~\\
		u_5^* \left\{
		\begin{array}{lllll}
			= u_{max} &\text{if}& \dis\frac{(\alpha_1-\alpha_2)S}{N}(\zeta_{ia,s}I_a + \zeta_{is,s}I_s + \zeta_{h,s}H) + (\alpha_6-\alpha_1)\tau_{r,s}R > -b_5(S+R)e^{-\sigma t} \\
			\in (0, u_{max}) &\text{if}& \dis\frac{(\alpha_1-\alpha_2)S}{N}(\zeta_{ia,s}I_a + \zeta_{is,s}I_s + \zeta_{h,s}H) + (\alpha_6-\alpha_1)\tau_{r,s}R = -b_5(S+R)e^{-\sigma t}\\
			= 0 &\text{if}& \dis\frac{(\alpha_1-\alpha_2)S}{N}(\zeta_{ia,s}I_a + \zeta_{is,s}I_s + \zeta_{h,s}H) + (\alpha_6-\alpha_1)\tau_{r,s}R < -b_5(S+R)e^{-\sigma t}
		\end{array}
		\right.\\~\\
		u_6^* \left\{
		\begin{array}{lllll}
			= u_{max} &\text{if}& (\alpha_2-\alpha_3)\tau_{e,ia}E + [(\alpha_3-\alpha_4)\tau_{ia,is}+(\alpha_3-\alpha_6)\tau_{ia,r}]I_a > -b_6(E+I_a)e^{-\sigma t} \\
			\in (0, u_{max}) &\text{if}& (\alpha_2-\alpha_3)\tau_{e,ia}E + [(\alpha_3-\alpha_4)\tau_{ia,is}+(\alpha_3-\alpha_6)\tau_{ia,r}]I_a = -b_6(E+I_a)e^{-\sigma t}\\
			= 0 &\text{if}& (\alpha_2-\alpha_3)\tau_{e,ia}E + [(\alpha_3-\alpha_4)\tau_{ia,is}+(\alpha_3-\alpha_6)\tau_{ia,r}]I_a < -b_6(E+I_a)e^{-\sigma t}
		\end{array}
		\right.\\~\\
		u_7^* \left\{
		\begin{array}{lllll}
			= u_{max} &\text{if}& [(\alpha_4-\alpha_5)\tau_{is,h} + (\alpha_4-\alpha_6)\tau_{is,r}+(\alpha_4-\alpha_7)\tau_{is,d}]I_s > -b_7I_se^{-\sigma t} \\
			\in (0, u_{max}) &\text{if}& [(\alpha_4-\alpha_5)\tau_{is,h} + (\alpha_4-\alpha_6)\tau_{is,r}+(\alpha_4-\alpha_7)\tau_{is,d}]I_s = -b_7I_se^{-\sigma t}\\
			= 0 &\text{if}& [(\alpha_4-\alpha_5)\tau_{is,h} + (\alpha_4-\alpha_6)\tau_{is,r}+(\alpha_4-\alpha_7)\tau_{is,d}]I_s < -b_7I_se^{-\sigma t}
		\end{array}
		\right.\\~\\
		u_8^* \left\{
		\begin{array}{lllll}
			= u_{max} &\text{if}& [(\alpha_4-\alpha_5)\tau_{is,h} + (\alpha_4-\alpha_6)\tau_{is,r}+(\alpha_4-\alpha_7)\tau_{is,d}]I_s > -b_8I_se^{-\sigma t} \\
			\in (0, u_{max}) &\text{if}& [(\alpha_4-\alpha_5)\tau_{is,h} + (\alpha_4-\alpha_6)\tau_{is,r}+(\alpha_4-\alpha_7)\tau_{is,d}]I_s = -b_8I_se^{-\sigma t}\\
			= 0 &\text{if}& [(\alpha_4-\alpha_5)\tau_{is,h} + (\alpha_4-\alpha_6)\tau_{is,r}+(\alpha_4-\alpha_7)\tau_{is,d}]I_s < -b_8I_se^{-\sigma t}
		\end{array}
		\right.
\end{eqnarray*}
\begin{eqnarray*}
		u_9^* \left\{
		\begin{array}{lllll}
			= u_{max} &\text{if}& [(\alpha_6-\alpha_5)\tau_{h,r} + (\alpha_5-\alpha_7)\tau_{h,d}]H > -b_9He^{-\sigma t} \\
			\in (0, u_{max}) &\text{if}& [(\alpha_6-\alpha_5)\tau_{h,r} + (\alpha_5-\alpha_7)\tau_{h,d}]H = -b_9He^{-\sigma t}\\
			= 0 &\text{if}& [(\alpha_6-\alpha_5)\tau_{h,r} + (\alpha_5-\alpha_7)\tau_{h,d}]H < -b_9He^{-\sigma t}
		\end{array}
		\right.\\
	\end{eqnarray*}
	
	\noindent
	But we know that the optimal policy is achieved when the marginal benefit and marginal cost are equal. Hence, from the above characterization of the optimal control $u^*$ we deduce the effectiveness conditions \eqref{effectivenessCond}. This ends the proof of Theorem \ref{characTheoContr2}. In application, greater effectiveness will lead to faster epidemic control, thereby lowering the overall disease impact on the society.
\end{proof}

\subsection{Feasibility Objective Function}
This section is devoted to the practical considerations for implementing interventions, such as availability of healthcare infrastructure, ease of compliance by the population, and adaptability of measures. This will help evaluate how realistic it is to deploy our control measures. Like the effectiveness, we aim to maximize feasibility to ensure that interventions are not only theoretically effective but also realistically achievable in practice. Let's consider the feasibility objective function
\begin{equation}\label{objective3}
	J_3(u_1, u_2, \cdots, u_9) = \dis\int_{0}^{T} \left[ \lambda_1E + \lambda_2I_a + \lambda_3I_s + \lambda_4H + \dis\frac{1}{2}\dis\sum_{i=1}^{9}b_iu_i^2(t)\right]e^{-\sigma t}dt,
\end{equation}
where $\left\{\lambda_i\right\}_{i=1}^{4}$ are the balancing weight on the exposed, the asymptomatically infected, the symptomatically infected and the hospitalized individuals; $\left\{b_i\right\}_{i=1}^{9}$ are the balancing cost factors on the different controls $\left\{u_i\right\}_{i=1}^{9}$; $\sigma$ is a scaling factor and $T$ is the final time for control implementation.

\noindent
Our goal is to find an optimal control $u^* = (u_1^*, u_2^*, \cdots, u_9^*)$ such that 
\begin{equation}\label{optimiz3}
	J_3(u^*) = \max_{u\in\U}{J_3(u)}
\end{equation}
where $u=(u_1,u_2, \cdots, u_9) \in \U$.

\noindent
Analogous to Theorem \ref{control1}, we can establish that there exists an optimal control $u^*$ solution of Problem \eqref{optimiz3}.\\

\noindent
It now remains to characterize the optimal control $u^*$, solution of \eqref{optimiz3}. As in the previous cases, we will use Pontryagin's Maximum Principle to proceed.\\

\noindent
Let $\alpha_i, ~ i=1, \cdots,7$ be the adjoint variables corresponding to the state variables $S, E, I_a, I_s, H, R$ and $D$ respectively. The Hamiltonian $\H$ associated with our problem is defined by:

\begin{equation}\label{hamilton1}
	\begin{array}{rlllll}
		\H &=& \alpha_1\dis\frac{dS}{dt} + \alpha_2\dis\frac{dE}{dt} + \alpha_3\dis\frac{dI_a}{dt} + \alpha_4\dis\frac{dI_s}{dt} + \alpha_5\dis\frac{dH}{dt} + \alpha_6\dis\frac{dR}{dt} + \alpha_7\dis\frac{dD}{dt} \\
		&& + \left(\lambda_1E + \lambda_2I_a + \lambda_3 I_s + \lambda_4 H + \dis\frac{1}{2}\dis\sum_{i=1}^{9}b_iu_i^2\right)e^{-\sigma t},
	\end{array}
\end{equation}
and we have the following result.

\begin{theo}\label{characTheoContr3}
	The optimal control $u^*=(u_1^*,u_2^*, \cdots, u_9^*)$ solution of the optimization problem
	\begin{equation*}
		J_3(u^*) = \max_{u\in\U}{J_3(u)}
	\end{equation*}
	is characterized by:
	\begin{equation}\label{characSyst}
		u_i^* = \left\{
		\begin{array}{clllll}
			u_{i_{max}} &\text{if}& \omega_i^*<u_{i_{max}}\\
			\omega_i^* &\text{if}& u_{i_{max}} \leq \omega_i^*
		\end{array}
		\right., ~~ i=1, \cdots, 9
	\end{equation}
	where:
	$$\begin{array}{rllllll}
		\omega_1^* &=& \dis\frac{S^*(\alpha_2-\alpha_1)(\zeta_{ia,s}I_a^* + \zeta_{is,s}I_s^* + \zeta_{h,s}H^*)}{b_1N} \\~\\
		
		\omega_2^* &=& \dis\frac{S^*(\alpha_1-\alpha_2)(\zeta_{ia,s}I_a^* + \zeta_{is,s}I_s^* + \zeta_{h,s}H^*)}{b_2N}
	\end{array}$$
	
	$$\begin{array}{rllllll}
		\omega_3^* &=& \dis\frac{S^*(\alpha_2-\alpha_1)(\zeta_{ia,s}I_a^* + \zeta_{is,s}I_s^* + \zeta_{h,s}H^*)}{b_3N} + \dis\frac{(\alpha_1-\alpha_6)\tau_{r,s}R^*}{b_3}\\~\\
		
		\omega_4^* &=& \dis\frac{S^*(\alpha_2-\alpha_1)(\zeta_{ia,s}I_a^* + \zeta_{is,s}I_s^* + \zeta_{h,s}H^*)}{b_4N} + \dis\frac{(\alpha_1-\alpha_6)\tau_{r,s}R^*}{b_4}\\~\\
		
		\omega_5^* &=& \dis\frac{S^*(\alpha_2-\alpha_1)(\zeta_{ia,s}I_a^* + \zeta_{is,s}I_s^* + \zeta_{h,s}H^*)}{b_5N} + \dis\frac{(\alpha_1-\alpha_6)\tau_{r,s}R^*}{b_5}\\~\\
		
		\omega_6^* &=& \dis\frac{(\alpha_3-\alpha_2)\tau_{e,ia}E^* + (\alpha_4-\alpha_3)\tau_{ia,is}I_a^* + (\alpha_6-\alpha_3)\tau_{ia,r}I_a^*}{b_6}\\~\\
		
		\omega_7^* &=& \dis\frac{\left[(\alpha_5\tau_{is,h}+\alpha_6\tau_{is,r}+\alpha_7\tau_{is,d}) - \alpha_4(\tau_{is,r}+\tau_{is,h}+\tau_{is,d})\right]I_s^*}{b_7}
	\end{array}$$
	
	$$\begin{array}{rllllll}
		\omega_8^* &=& \dis\frac{\left[(\alpha_5\tau_{is,h}+\alpha_6\tau_{is,r}+\alpha_7\tau_{is,d}) - \alpha_4(\tau_{is,r}+\tau_{is,h}+\tau_{is,d})\right]I_s^*}{b_8}\\~\\
		
		\omega_9^* &=& \dis\frac{\left[(\alpha_5-\alpha_6)\tau_{h,r} + (\alpha_7-\alpha_5)\tau_{h,d}\right]H^*}{b_9}
	\end{array}$$
	with the septuple $(S^*, E^*, I_a^*, I_s^*, H^*, R^*, D^*)$ solution of the control system \eqref{contsyst1} and the adjoint variables $(\alpha_1, \alpha_2, \cdots, \alpha_7)$ satisfy the following adjoint system, with transversality conditions $\alpha_i(T) = 0, ~~ i=1,\cdots,7$
	
	\begin{equation}\label{adjointSyst1}
		\left\{
		\begin{array}{lllllll}
			\dis\frac{d\alpha_1}{dt} &=& \dis\frac{(\alpha_1-\alpha_2)}{N}(1-u_1+u_2-u_3-u_4-u_5)(\zeta_{ia,s}I_a + \zeta_{is,s}I_s + \zeta_{h,s}H) \\~\\		
			\dis\frac{d\alpha_2}{dt} &=& (\alpha_2-\alpha_3)(1-u_6)\tau_{e,ia} -\lambda_1e^{-\sigma t}\\~\\
			\dis\frac{d\alpha_3}{dt} &=& \dis\frac{(\alpha_1-\alpha_2)S}{N}(1-u_1+u_2-u_3-u_4-u_5)\zeta_{ia,s} + (1-u_6)\left[(\alpha_3-\alpha_4)\tau_{ia,is} \right.\\ 
			&& \left. + (\alpha_3-\alpha_6)\tau_{ia,r}\right] - \lambda_2e^{-\sigma t} \\~\\
			\dis\frac{d\alpha_4}{dt} &=& \dis\frac{(\alpha_1-\alpha_2)S}{N}(1-u_1+u_2-u_3-u_4-u_5)\zeta_{is,s} + (1-u_7-u_8)\left[(\alpha_4-\alpha_6)\tau_{is,r} \right.\\ 
			&& \left. + (\alpha_4-\alpha_5)\tau_{is,h} + (\alpha_4-\alpha_7)\tau_{is,d}\right] - \lambda_3e^{-\sigma t}\\~\\
			\dis\frac{d\alpha_5}{dt} &=& \dis\frac{(\alpha_1-\alpha_2)S}{N}(1-u_1+u_2-u_3-u_4-u_5)\zeta_{h,s} + (1-u_9)(\alpha_5-\alpha_7)\tau_{h,d} \\ 
			&& + u_9(\alpha_5-\alpha_6)\tau_{h,r} - \lambda_4e^{-\sigma t}\\~\\
			\dis\frac{d\alpha_6}{dt} &=& \dis(\alpha_6-\alpha_1)(1-u_3-u_4-u_5)\tau_{r,s}\\~\\
			\dis\frac{d\alpha_7}{dt} &=& 0.
		\end{array}
		\right.
	\end{equation}
\end{theo}

\noindent
\begin{proof}
	Consider the following Hamilton's equations, derived from the total energy of the system (called Hamiltonian):\\
	
	~\noindent
	$\dis\frac{\partial \H}{\partial S} = -\dis\frac{d\alpha_1}{dt},
	\dis\frac{\partial \H}{\partial E} = -\dis\frac{d\alpha_2}{dt},
	\dis\frac{\partial \H}{\partial I_a} = -\dis\frac{d\alpha_3}{dt},
	\dis\frac{\partial \H}{\partial I_s} = -\dis\frac{d\alpha_4}{dt},
	\dis\frac{\partial \H}{\partial H} = -\dis\frac{d\alpha_5}{dt},
	\dis\frac{\partial \H}{\partial R} = -\dis\frac{d\alpha_6}{dt},
	\dis\frac{\partial \H}{\partial D} = -\dis\frac{d\alpha_7}{dt}
	$.\\
	
	\noindent
	The expanded form of the Hamiltonian is given below:
	\begin{eqnarray*}
		\H &=& \alpha_1 \left[ -\left(1-u_1+u_2-u_3-u_4-u_5\right) \left(\zeta_{ia,s}I_a + \zeta_{is,s}I_s + \zeta_{h,s}H\right)\dis\frac{S}{N} + (1-u_3-u_4-u_5)\tau_{r,s}R \right] \\
		&& + \alpha_2 \left[ \left(1-u_1+u_2-u_3-u_4-u_5\right) \left(\zeta_{ia,s}I_a + \zeta_{is,s}I_s + \zeta_{h,s}H\right)\dis\frac{S}{N} - (1-u_6)\tau_{e,ia}E \right] \\
		&& + \alpha_3 \left[ (1-u_6)\tau_{e,ia}E - (1-u_6) (\tau_{ia,is} + \tau_{ia,r})I_a \right]\\
		&& + \alpha_4 \left[ (1-u_6)\tau_{ia,is}I_a - (1-u_7-u_8)(\tau_{is,r} + \tau_{is,h} + \tau_{is,d})I_s \right] \\
		&& + \alpha_5 \left[ (1-u_7-u_8)\tau_{is,h}I_s - (1-u_9)\tau_{h,d}H -u_9\tau_{h,r}H\right] \\
		&& + \alpha_6 \left[ (1-u_7-u_8)\tau_{is,r}I_s + (1-u_6)\tau_{ia,r}I_a + u_9\tau_{h,r}H \right. \\
		&& \left. -(1-u_3-u_4-u_5)\tau_{r,s}R \right] + \alpha_7 \left[ (1-u_9)\tau_{h,d}H + (1-u_7-u_8)\tau_{is,d}I_s \right]\\
		&& + \left[\lambda_1E + \lambda_2I_a + \lambda_3 I_s + \lambda_4 H + \dis\frac{1}{2}\dis\sum_{i=1}^{9}b_iu_i^2\right]e^{-\sigma t}.
	\end{eqnarray*}
	
	\noindent
	The adjoint system is derived from the Hamilton's equations and given as follows:
	
	\begin{eqnarray*}
		\dis\frac{d \alpha_1}{dt} &=& \dis\frac{(\alpha_1-\alpha_2)}{N}(1-u_1+u_2-u_3-u_4-u_5)(\zeta_{ia,s}I_a + \zeta_{is,s}I_s + \zeta_{h,s}H)\\
		\dis\frac{d \alpha_2}{dt} &=& (\alpha_2-\alpha_3)(1-u_6)\tau_{e,ia} - \lambda_1e^{-\sigma t} \\
		\dis\frac{d \alpha_3}{dt} &=& \dis\frac{(\alpha_1-\alpha_2)S}{N}(1-u_1+u_2-u_3-u_4-u_5)\zeta_{ia,s} + (1-u_6)\left[(\alpha_3-\alpha_4)\tau_{ia,is}\right.\\
		&& \left. + (\alpha_3-\alpha_6)\tau_{ia,r}\right] - \lambda_2e^{-\sigma t}\\
		\dis\frac{d \alpha_4}{dt} &=& \dis\frac{(\alpha_1-\alpha_2)S}{N}(1-u_1+u_2-u_3-u_4-u_5)\zeta_{is,s} + (1-u_7-u_8)\left[(\alpha_4-\alpha_6)\tau_{is,r}\right.\\
		&& \left. + (\alpha_4-\alpha_5)\tau_{is,h} + (\alpha_4-\alpha_7)\tau_{is,d}\right] - \lambda_3e^{-\sigma t}\\
		\dis\frac{d \alpha_5}{dt} &=& \dis\frac{(\alpha_1-\alpha_2)S}{N}(1-u_1+u_2-u_3-u_4-u_5)\zeta_{h,s} + (1-u_9)(\alpha_5-\alpha_7)\tau_{h,d}\\
		&& + u_9(\alpha_5-\alpha_6)\tau_{h,r} - \lambda_4e^{-\sigma t} \\
		\dis\frac{d \alpha_6}{dt} &=& (\alpha_6-\alpha_1)(1-u_3-u_4-u_5)\tau_{r,s}\\
		\dis\frac{d \alpha_7}{dt} &=& 0.
	\end{eqnarray*}
	
	\noindent
	According to Remark \ref{Tranversality}, we have: $\alpha_i(T) = 0, ~~ i=1,\cdots, 7$.\\
	
	\noindent
	Now, for each $u_i,~~ i=1,\cdots,9$, we need to solve the following
	$$\dis\frac{\partial \H}{\partial u_i} = 0$$
	to complete the characterization of the optimal control.
	
	\noindent
	We have:\\
	$$\dis\frac{\partial \H}{\partial u_1} = (\alpha_1-\alpha_2)(\zeta_{ia,s}I_a + \zeta_{is,s}I_s + \zeta_{h,s}H)\dis\frac{S}{N} + b_1u_1e^{-\sigma t}.$$
	
	\noindent
	Solving $\dis\frac{\partial \H}{\partial u_1} = 0$ for $u_1$, we obtain:

	$$\tilde{u}_1 = \dis\frac{(\alpha_2-\alpha_1)(\zeta_{ia,s}I_a + \zeta_{is,s}I_s + \zeta_{h,s}H)Se^{\sigma t}}{b_1N}.$$
	
	\noindent
	Taking into account the admissible control set, if we set $u_{1_{max}} = \dis\max\{u_1\}$ for all $u_1\in [0,u_{max}]$, then it is clear that the optimal value of the control variable $u_1$ is given by:
	$$u_1^* = \max \left\{\omega_1^*, u_{1_{max}}\right\},$$
	where
	$$\omega_1^* = \dis\frac{S^*(\alpha_2-\alpha_1)(\zeta_{ia,s}I_a^* + \zeta_{is,s}I_s^* + \zeta_{h,s}H^*)}{b_1N}$$
	Hence:
	
	$$u_1^* = \left\{
	\begin{array}{rlll}
		u_{1_{max}} &\text{if}& \omega_1^* < u_{1_{max}} \\
		\omega_1^* &\text{if}& u_{1_{max}} \leq \omega_1^*
	\end{array}
	\right.$$
	Similarly, we have:
	
	$$u_i^* = \left\{
	\begin{array}{rlll}
		u_{i_{max}} &\text{if}& \omega_i^* < u_{i_{max}} \\
		\omega_i^* &\text{if}& u_{i_{max}} \leq \omega_i^*
	\end{array}
	\right., ~~ i=2,\cdots, 9$$
	where:
	
	\begin{eqnarray*}
		\omega_2^* &=& \dis\frac{S^*(\alpha_1-\alpha_2)(\zeta_{ia,s}I_a^* + \zeta_{is,s}I_s^* + \zeta_{h,s}H^*)}{b_2N}\\~\\
		\omega_3^* &=& \dis\frac{S^*(\alpha_2-\alpha_1)(\zeta_{ia,s}I_a^* + \zeta_{is,s}I_s^* + \zeta_{h,s}H^*)}{b_3N} + \dis\frac{(\alpha_1-\alpha_6)\tau_{r,s}R^*}{b_3}\\~\\
		\omega_4^* &=& \dis\frac{S^*(\alpha_2-\alpha_1)(\zeta_{ia,s}I_a^* + \zeta_{is,s}I_s^* + \zeta_{h,s}H^*)}{b_4N} + \dis\frac{(\alpha_1-\alpha_6)\tau_{r,s}R^*}{b_4}\\~\\
		\omega_5^* &=& \dis\frac{S^*(\alpha_2-\alpha_1)(\zeta_{ia,s}I_a^* + \zeta_{is,s}I_s^* + \zeta_{h,s}H^*)}{b_5N} + \dis\frac{(\alpha_1-\alpha_6)\tau_{r,s}R^*}{b_5}\\~\\
		\omega_6^* &=& \dis\frac{(\alpha_3-\alpha_2)\tau_{e,ia}E^* + (\alpha_4-\alpha_3)\tau_{ia,is}I_a^* + (\alpha_6-\alpha_3)\tau_{ia,r}I_a^*}{b_6}
	\end{eqnarray*}
	
	\begin{eqnarray*}
		\omega_7^* &=& \dis\frac{\left[(\alpha_5\tau_{is,h} + \alpha_6\tau_{is,r} + \alpha_7\tau_{is,d}) - \alpha_4(\tau_{is,r} + \tau_{is,h} + \tau_{is,d})\right]I_s^*}{b_7} \\~\\
		\omega_8^* &=& \dis\frac{\left[(\alpha_5\tau_{is,h} + \alpha_6\tau_{is,r} + \alpha_7\tau_{is,d}) - \alpha_4(\tau_{is,r} + \tau_{is,h} + \tau_{is,d})\right]I_s^*}{b_8} \\~\\
		\omega_9^* &=& \dis\frac{\left[(\alpha_5-\alpha_6)\tau_{h,r} + (\alpha_7 - \alpha_5)\tau_{h,d}\right]H^*}{b_9}.
	\end{eqnarray*}
	
	\noindent
	This ends the proof of Theorem \ref{characTheoContr3}.
\end{proof}


\section{Concluding remarks}\label{conclusion}
This study contributes to the theoretical advancement of Pontryagin's Maximum Principle in epidemiological modeling by providing a structured evaluation of single-criterion interventions. In this article, we explored the optimization of COVID-19 control strategies using Pontryagin's Maximum Principle, focusing on individual criteria: cost, effectiveness, and feasibility. This approach offers valuable insights into how these factors could influence public health outcomes. For instance, the cost criterion evaluates financial sustainability, ensuring that interventions remain economically viable. The effectiveness criterion targets minimizing infections and deaths to control the overall impact of the disease. The feasibility criterion emphasizes practical implementation, considering real-world conditions and compliance variability. Analyzing these criteria independently allowed us to evaluate their contributions, offering a clearer understanding of the trade-offs inherent in single-objective optimization.\\

\noindent
Despite these insights, single-criterion strategies might have inherent limitations: the cost-focused approach may lack sufficient impact, highly effective interventions may be impractical for widespread implementation, and the feasible approach may fall short in reducing the socioeconomic impact of the disease. As a next step, we then plan to conduct numerical simulations to validate and further refine the proposed model. These simulations will provide a deeper understanding of how the theoretical findings translate into practical outcomes for various scenarios. By incorporating real-world data, the simulations will evaluate the optimality of interventions for each individual criterion. They will allow us to explore the potential of multiple objectives in a more comprehensive manner. A natural progression of this work then involves multi-objective optimization, where cost, effectiveness, and feasibility are addressed simultaneously to design more robust and sustainable intervention strategies. Such an approach could better capture the complexity of public health decision-making, fostering comprehensive solutions to pandemic challenges.\\

\section*{Statements and Declarations}
Data Availability: No data were generated or analyzed during this study.

\noindent The authors declare that there is no conflict of interest regarding the publication of this article.

\section*{Competing Interests and Funding}
The authors declare that they have no competing interests.
This work was funded by the U.S. Department of Energy (DOE) under grant DE-AC02-05CH11231 through the Biopreparedness Research Virtual Environment (BRaVE) program,  "EMERGE: ExaEpi for Elucidating Multiscale Ecosystem Complexities for Robust, Generalized Epidemiology."



\end{document}